\documentclass[a4paper,11pt]{amsart}
      
\usepackage{amsmath,amsfonts,amssymb}
\usepackage{enumerate} 
\usepackage{amsthm}
\def\bound{\widetilde{M}_\para^\e}
\def\GO{\mathcal{O}}

\def\ga{{g}}
\def\C{{\mathbb C}}
\def\R{{\mathbb R}}
\def\N{{\mathbb N}}

\def\T{{\mathbb T}}

\def\virgp{\raise 2pt\hbox{,}}
\def\bv{{\bf v}}
\def\bu{{\bf u}}

\def\({\left(}
\def\){\right)}
\def\<{\left\langle}
\def\>{\right\rangle}
\def\le{\leqslant}
\def\ge{\geqslant}

\def\Tend#1#2{\mathop{\longrightarrow}\limits_{#1\rightarrow#2}}

\def\d{{\partial}}

\def\e{\varepsilon}

\def\F{\mathcal F}

\def\O{\mathcal O}
\def\eik{\phi_{\ei}}
\def\ext{V_{\rm ext}}
\def\poi{V_{\rm p}}

\DeclareMathOperator{\RE}{Re}
\DeclareMathOperator{\IM}{Im}
\DeclareMathOperator{\cn}{\mathrm{div}}

\DeclareMathOperator{\ei}{eik}

\def\defn{\mathrel{:=}}

\def\id{I}
\def\la{\left\lvert}

\def\le{\leqslant}
\def\les{\lesssim}

\def\para{h}

\def\ra{\right\rvert}

\def\scal#1#2{\langle \, {#1} \hspace{1pt},\hspace{1pt} {#2} \,\rangle}

\theoremstyle{plain}
\newtheorem{theorem}{Theorem}[section]
\newtheorem{lemma}[theorem]{Lemma}
\newtheorem{corollary}[theorem]{Corollary}
\newtheorem{proposition}[theorem]{Proposition}
\newtheorem{hyp}{Assumption}

\theoremstyle{definition}

\newtheorem*{notation}{Notation}

\theoremstyle{remark}
\newtheorem{remark}[theorem]{Remark}
\newtheorem*{remark*}{Remark}
\newtheorem*{example}{Example}

\numberwithin{equation}{section}

\begin{document}

\title[Semi-classical limit of Schr\"odinger--Poisson]{Semi-classical
 limit of Schr\"odinger--Poisson equations in space dimension $n\ge 3$}

\author[T. Alazard]{Thomas Alazard}
\address{Universit\'e Paris-Sud\\ Math\'ematiques\\ B\^at. 425\\ 91405
  Orsay cedex\\ France}
\email{Thomas.Alazard@math.u-psud.fr}
\author[R. Carles]{R{\'e}mi Carles}
\address{Wolfgang Pauli Institute, c/o Inst.~f.~Math., Universit\"at Wien,
        Nordbergstr.~15, A-1090 Wien, Austria\footnote{On leave from
  MAB, UMR 5466 CNRS, Bordeaux, France.}}
\email{Remi.Carles@math.cnrs.fr}
\thanks{Supports by European network HYKE,
 funded  by the EC as contract HPRN-CT-2002-00282, by Centro de
 Matem\'atica e Aplica\c c\~oes Fundamentais (Lisbon), funded by
 FCT as contract POCTI-ISFL-1-209, and by the ANR project SCASEN,
 are acknowledged.} 

\begin{abstract}
We prove the existence of solutions to the Schr\"odinger--Poisson
system on a time interval independent of the Planck constant, when the
doping profile does not necessarily decrease at infinity, in the
presence of a subquadratic external potential. The lack of
integrability of the doping profile is resolved by working in Zhidkov
spaces, in space dimension at least three. We infer that the main
quadratic quantities (position density and modified momentum density)
converge strongly  as the Planck constant goes to zero. When the
doping profile is integrable, we prove pointwise convergence.  
\end{abstract}

\maketitle

\section{Introduction}\label{sec:intro}

We consider the semi-classical limit $\e\to 0$ of the
Schr\"odinger--Poisson system: 
\begin{align}
      i\e \d_t u^\e +\frac{\e^2}{2}\Delta u^\e &= \ext u^\e +
      \poi^\e u^\e ,\quad (t,x)\in \R\times \R^n,\label{eq:schrod}\\
\Delta \poi^\e &= q\(|u^\e|^2-c\), 
\label{eq:poisson}\\
u^\e_{\mid t=0} &= a_0^\e(x)e^{i\Phi_0(x)/\e},\label{eq:CI}
\end{align}
where $\ext=\ext(t,x)$ is an external potential (harmonic
potential for instance),  $c=c(x)$ is a \emph{doping profile} (or
\emph{impurity},  
\emph{background ions}), and $q\in \R$ represents an electric charge;
$\ext$, $c$ and $q$ are data of the problem (see
e.g. \cite{MaRiSc}). We consider the case where the space dimension is
$n\ge 3$. This is due to a lack of control of low frequencies for the
Poisson equation \eqref{eq:poisson} when $n\le 2$. 
\smallbreak

The conditions we impose to solve the Poisson equation
\eqref{eq:poisson} will be given according to the different cases we
consider. 
\smallbreak

The doping profile $c$ is supposed to be bounded, and does not
necessarily goes to zero at infinity (see
Assumption~\ref{hyp:lin} or Assumption~\ref{hyp:general}
below). Suppose for instance 
that $c\equiv 1$. Then \eqref{eq:schrod}--\eqref{eq:poisson} is
reminiscent of the 
Gross-Pitaevskii equation (see e.g. \cite{LinZhang,PG05} and
references therein):
\begin{equation}
  \label{eq:gp}
  i\e \d_t u^\e +\frac{\e^2}{2}\Delta u^\e = \( |u^\e|^2-1\)u^\e.
\end{equation}
For this equation, the Hamiltonian structure yields, at least
formally:
\begin{equation*}
  \frac{d}{dt}\(\|\e \nabla u^\e(t)\|_{L^2}^2 +
  \left\||u^\e(t)|^2-1\right\|_{L^2}^2  \)=0.
\end{equation*}
A natural
space to study the Cauchy problem associated to \eqref{eq:gp}  is
therefore the energy space 
\begin{equation*}
  E=\{ u \in H^1_{\rm loc}(\R^n)\ ;\ \nabla u\in L^2(\R^n), \ |u|^2-1
  \in L^2(\R^n)\}. 
\end{equation*}
For this quantity to be well defined, one cannot assume that $u^\e$ is
in $L^2(\R^n)$; morally, the modulus of $u^\e$
goes to one at infinity. To study solutions which are bounded, but not
in $L^2(\R^n)$, P.~E.~Zhidkov introduced in the
one-dimensional case in \cite{Zhidkov} (see also \cite{ZhidkovLNM}):
\begin{equation*}
  X^s(\R^n) = \{  u \in L^\infty(\R^n)\ ;\ \nabla u \in
  H^{s-1}(\R^n)\},\quad s>n/2. 
\end{equation*}
The study of these spaces was 
generalized in the multidimensional case by C.~Gallo
\cite{Gallo}. They make it possible to consider solutions to
\eqref{eq:gp} whose modulus has a non-zero limit as $|x|\to
\infty$, but not necessarily satisfying $|u^\e(t,\cdot)|^2 -1\in
L^2(\R^n)$. 

Recently, P.~G\'erard \cite{PG05} solved the Cauchy problem for the
Gross-Pitaevskii equation in the more natural space $E$, in space
dimensions two and three. The main
novelty consists in working 
with distances instead of norms, in order to apply a fixed point
argument in $E$. In particular, the constraint  $|u^\e(t,\cdot)|^2 -1\in
L^2(\R^n)$ is satisfied (and propagated). 

We have to face a similar issue, when solving
the Poisson equation. Mimicking the approach of \cite{LinZhang,PG05},
it is natural to work with the property:
\begin{equation*}
  |u^\e(t,\cdot)|^2-c(\cdot)\in L^2(\R^n).
\end{equation*}
We shall always assume that this holds at time $t=0$. 
We prove that this property holds on $[0,T]$ for some $T>0$
independent of $\e$, provided that we consider an external
potential whose unbounded part is linear in $x$. However, our analysis
shows that in the 
presence of a quadratic external potential, this property is not
relevant off $t=0$   (see Section~\ref{sec:eik}).
\smallbreak
 
Note that we make no assumption on the sign of $q$ (which models the
charge of the element considered in a semiconductor device). This is
in sharp contrast with the mathematical analysis of the semi-classical
limit of the  nonlinear
Schr\"odinger equation.  When the Poisson term $\poi^\e(t,x)u^\e$ 
is replaced with the nonlinear term $f(\left\lvert u^{\e}\right\rvert
^{2})u^\e$, E.~Grenier \cite{Grenier98} proposed a 
 strategy to obtain a phase/amplitude 
representation of the solution $u^\e$. This leads to study  
a quasi-linear system whose principal part writes:
$$
\square_{f'}\defn
\d_{t}^{2}-\cn\bigl(f'(\left\lvert u^\e\right\rvert
^2)\nabla\cdot\bigr). 
$$
Hence, to prove that the Cauchy problem is well-posed, one 
has to assume that the nonlinearity is defocusing and cubic at the
origin ($f'>0$), except for  
analytic initial data \cite{PGX93}, for which one can solve elliptic
evolution equations. Here, 
we are not restricted to the case when $q>0$. As will be clear below,
the reason is that the   
quasi-linear operator $\square_{f'}$ is replaced with the
\emph{semi-linear} operator  
$\d_{t}^{2}-q\Delta^{-1}\nabla((\la u^\e\ra^2-1)\cn \cdot)$. 

\begin{notation}
Recall that for $s>n/2$, Zhidkov spaces are defined by\footnote{For
  general $s>0$, another definition is used, see \cite{Gallo}.}:  
\begin{equation*}
  X^s(\R^n)=\{  u \in L^\infty(\R^n)\ ;\ \nabla u \in
  H^{s-1}(\R^n)\}\cdot
\end{equation*}
We denote
\begin{equation*}
  \|u\|_{X^s}\defn \|u\|_{L^\infty} + \|\nabla u\|_{H^{s-1}}.
\end{equation*}
We write 
$H^s=H^s(\R^n)$, $X^s=X^s(\R^n)$, 
$H^\infty \defn \cap_{s\in \N}H^s$, $X^\infty \defn \cap_{s\in \N}X^s$. 
We do not use specific notations for vector-valued functions: for instance, 
we write abusively $\nabla^2 f\in H^\infty$ when 
$\d^2_{jk} f\in H^\infty$ for every $1\le j,k\le n$.
\end{notation}
\begin{remark}
Zhidkov spaces contain  all the functions of the form
\begin{equation*}
\gamma+v,\text{ with }\gamma={\rm Const.}\in \C \text{ and }v\in H^s(\R^n).
\end{equation*}
The converse is not true, as shown by the following example:
\begin{equation*}
u(x)=\frac{x_{1}}{1+|x|^2}\virgp \quad x=(x_{1},x_{2},x_{3})\in \R^3.
\end{equation*}
On the other hand, if $n\ge 3$ and $u\in X^s$ for some $s>n/2$, then
there exists $\gamma\in \C$ such that $u-\gamma \in
L^{\frac{2n}{n-2}}(\R^n)$ (see Lemma~\ref{lem:hormander} below).  
\end{remark}

In this paper, we consider the system \eqref{eq:schrod}--\eqref{eq:CI}
in three cases:
\begin{itemize}
\item The external potential
  and the initial phase are sub-linear in $x$, and the mobility $c$ is
  in Zhidkov spaces (Part~\ref{part:1}). 
\item The external potential
  and the initial phase are sub-quadratic in $x$, and $c$ is
  a short range perturbation of a non-zero constant (Part~\ref{part:2}).
\item The mobility is integrable, and the external potential
  and the initial phase are sub-quadratic in $x$ (Part~\ref{part:3}).
\end{itemize}
In the first two cases, we construct a solution to
\eqref{eq:schrod}--\eqref{eq:CI} in Zhidkov spaces, and describe the
asymptotic behavior of the main quadratic observables as $\e\to 0$. In
the last case, we construct a solution in Sobolev spaces, and give
pointwise asymptotics of the solution as $\e\to 0$. 
\smallbreak

In this introduction, we describe more precisely the results
corresponding to the first case. We emphasize the fact that if we
simply assume $V_{\rm 
  ext}\in C(\R;H^\infty)$ and $\Phi_0\in H^\infty$, then our analysis
becomes much simpler. The unboundedness of $V_{\rm
  ext}$ and $\Phi_0$ require some geometrical description that
complicates the technical approach. Yet, this makes our assumptions
more physically relevant (see e.g. \cite{GosseMauser} and references
therein).  
\begin{hyp}\label{hyp:lin} Recall that $n\ge 3$.\\
\noindent$\bullet$  \emph{External potential:}  $\ext\in
C^\infty(\R\times\R^n)$ writes 
\begin{equation*}
  \ext (t,x)= E(t)\cdot x + V_{\rm pert}(t,x), \text{ with }E\in
  C^\infty(\R) \text{ and }\nabla V_{\rm pert}\in C(\R;H^\infty).
\end{equation*}

\noindent$\bullet$  \emph{Doping profile:}  $c\in
X^\infty$.\\

\noindent$\bullet$ \emph{Initial amplitude:} $a_0^\e(x) =a_0(x) +
r^\e(x)$, 
where $ a_0\in X^\infty$ is such that $|a_0|^2- c\in
L^2(\R^n)$, and $r^\e\in H^\infty$, with 
\begin{equation*}
  \|r^\e\|_{H^s}\Tend \e 0 0,\quad \forall s\ge 0.
\end{equation*}

\noindent$\bullet$ \emph{Initial phase:} we have $\Phi_0\in
C^\infty(\R^n)$ with
\begin{equation*}
  \Phi_0(x) = \alpha_0\cdot x +\phi_0(x),\text{ with }\alpha_0 \in
  \R^n\text{ and }\nabla \phi_0\in H^\infty.
\end{equation*}
\end{hyp}
\begin{lemma}\label{lem:eik0}
  Under the Assumption~\ref{hyp:lin}, there exists a unique
  solution $\eik\in C^\infty(\R\times\R^n)$ to:
  \begin{equation}
    \label{eq:eik0}
    \d_t \eik +\frac{1}{2}|\nabla \eik|^2 +E(t)\cdot x=0\quad ;\quad
    \eik(0,x)=\alpha_0\cdot x+\beta_0\, .
  \end{equation}
This solution is given explicitly by $\eik(t,x)= \alpha(t)\cdot x
+\beta(t)$, where:
\begin{equation*}
  \alpha(t)=\alpha_0 -\int_0^t E(\tau)d\tau\quad ;\quad \beta(t)=
  \beta_0-\frac{1}{2}\int_0^t \alpha(\tau)^2d\tau.
\end{equation*}
\end{lemma}
We skip the proof of this lemma; a more general result is proved
in Section~\ref{sec:eik}. We will see that if $V_{\rm ext}$ and/or
$\Phi_0$ have a quadratic dependence on $x$, then we have to consider
an eikonal phase $\eik$ which is quadratic in $x$. 
\begin{theorem}\label{theo:existence}
Let Assumption~\ref{hyp:lin} be  satisfied. 
 There exists
$T>0$ independent 
  of $\e \in ]0,1]$ and a  solution $u^\e \in
  L^\infty([0,T]\times \R^n)$
 to   \eqref{eq:schrod}-\eqref{eq:CI}, with
 \begin{equation*}
   \nabla\poi^\e(t,x)\to 0\ \text{ as
      }|x|\to \infty,\quad \poi^\e(t,0)=0,
 \end{equation*}
and such that $|u^\e|^2-c \in
  L^\infty([0,T];L^2)$. Moreover, one can write $u^\e = a^\e e^{i(\eik
    +\phi^\e)/\e}$, where:
  \begin{itemize}
  \item $a^\e \in  C^\infty([0,T]\times\R^n)\cap C([0,T];X^\infty)$,
  and $|a^\e|^2-c \in 
  C([0,T];L^2)$. 
  \item $\phi^\e \in C^\infty([0,T]\times\R^n)$ and $\nabla\phi^\e \in
  C([0,T];X^\infty)$. 
  \item We have the following uniform estimate: for every $s>n/2$,
  there exists $M_s$ \emph{independent of $\e \in ]0,1]$} such that
  \begin{equation*}
    \| a^\e\|_{L^\infty (0,T;X^s)}
    +\left\||a^\e|^2-c\right\|_{L^\infty (0,T;L^2)}+ \|\nabla
    \phi^\e\|_{L^\infty (0,T;X^s)}\le M_s. 
  \end{equation*}
  \end{itemize}
\end{theorem}
\begin{remark}
 We could not prove a uniqueness result for $u^\e$. 
 \end{remark}
 \begin{remark}
The above conditions to solve the Poisson equation are similar to
those given in \cite{ZZM}. We explain at the end of
Section~\ref{sec:cvschemelin} why in our framework, we cannot impose
$\poi^\e(t,x)\to 0$ as $|x|\to \infty$ (as in
\cite{BrezziMarkowich,ZhangSIMA} for instance). 
\end{remark}
Besides the uniform bounds, even the existence of such a solution
$u^\e$ is new. First, 
the presence of the external potential seems to have never been
studied rigorously before. As we already mentioned, this makes the
proof more technically involved. Next, in most
of the previous studies, $u^\e$ is supposed to be in $L^2$:  see e.g.
\cite{CastoM3AS97,Sulem}. In \cite{ZhangSIMA}, the author considers
the case $c\in L^1\cap H^s$. As we will see in
Section~\ref{sec:integrable}, this case makes the analysis easier, and also
makes it possible to have $u^\e \in L^2$. 
The main difficulty in the analysis lies 
in the fact that when $c\not \in L^1(\R^n)$, the condition
$|u^\e|^2-c\in L^2(\R^n)$ is somehow ``more nonlinear'', 
as in \cite{PG05}. 
\smallbreak
 
The general idea to prove
Theorem~\ref{theo:existence} 
consists in adapting the idea of \cite{Grenier98}: with techniques
from the hyperbolic theory, we construct a solution to 
\begin{equation}
  \label{eq:systcomplet}
\begin{aligned}
  \d_t \Phi^\e +\frac{1}{2}|\nabla \Phi^\e|^2 + \ext + \poi^\e&=0
  \quad ;\quad \Phi^\e_{\mid t=0}=\Phi_0.\\
\d_t a^\e +\nabla \Phi^\e \cdot \nabla a^\e +\frac{1}{2}a^\e \Delta
  \Phi^\e &= i\frac{\e}{2}\Delta a^\e\quad ;\quad a^\e_{\mid t=0}=a_0^\e.\\
\Delta \poi^\e =q\big(|a^\e|^2-c\big)\quad ;\quad &\nabla\poi^\e(t,x)\Tend
  {|x|}{\infty}0,\ \poi^\e(t,0)=0. 
\end{aligned}
\end{equation}
Following \cite{CaBKW}, we write $\Phi^\e =\eik + \phi^\e$:
with the unknown $(a^\e,\phi^\e)$, \eqref{eq:systcomplet} becomes (we
keep the term $\Delta \eik$ which is zero here, for future references):
\begin{equation}
  \label{eq:systcomplet2}
\begin{aligned}
  \d_t \phi^\e + \nabla \eik\cdot
  \nabla \phi^\e+\frac{1}{2}|\nabla \phi^\e|^2 +V_{\rm pert} + \poi^\e&=0
  \ ;\ \phi^\e_{\mid t=0}=\phi_0.\\
\d_t a^\e +\nabla \(\phi^\e+\eik\) \cdot \nabla a^\e +\frac{1}{2}a^\e \Delta
  \(\phi^\e+\eik\) &= i\frac{\e}{2}\Delta a^\e\ ;\ a^\e_{\mid
  t=0}=a_0^\e.\\ 
\Delta \poi^\e =q\big(|a^\e|^2-c\big)\quad ;\quad \nabla\poi^\e(t,x)\Tend
  {|x|}{\infty}0&,\ \poi^\e(t,0)=0. 
\end{aligned}
\end{equation}
Proving the existence and uniqueness of solution to
\eqref{eq:systcomplet2} as we do in the proof of
Theorem~\ref{theo:existence} is one of the main results of this
paper. Because of the difficulties pointed out above, and the fact that
one can easily be mistaken by using the usual approach, we give full
details for the construction of the solution to
\eqref{eq:systcomplet2}. 
Passing formally to the limit, it is natural to consider:
\begin{equation}
  \label{eq:systlim1}
 \begin{aligned}
  \d_t \phi + \nabla \eik\cdot
  \nabla \phi+\frac{1}{2}|\nabla \phi|^2 +V_{\rm pert} + \poi&=0
  \ ;\ \phi_{\mid t=0}=\phi_0.\\
\d_t a +\nabla \(\phi+\eik\) \cdot \nabla a +\frac{1}{2}a \Delta
  \(\phi+\eik\) &= 0\ ;\ a_{\mid
  t=0}=a_0.\\ 
\Delta \poi =q\big(|a|^2-c\big)\quad ;\quad \nabla\poi(t,x)\Tend
  {|x|}{\infty}0&,\ \poi(t,0)=0.
\end{aligned} 
\end{equation}

\begin{notation}
The symbol $\les$ stands for $\le$ up to a positive, multiplicative constant 
which depends only on parameters that are considered fixed. \\
We shall also denote $L^\infty_T Y$  for $L^\infty ([0,T];Y)$. 
\end{notation}
\begin{theorem}\label{theo:convergence}
  Under Assumption~\ref{hyp:lin}, there
  exists a smooth solution $(a,\phi)$ of \eqref{eq:systlim1} 
such that 
  $a,\nabla\phi \in C([0,T],X^\infty)$, $|a|^2-c \in
  C([0,T],L^2)$, 
  and 
  $$
  \|a^\e-a\|_{L^\infty_{T}H^s}+ \|\nabla(\phi^\e
  -\phi)\|_{L^\infty_{T}X^s}\Tend \e 0 0, \qquad \forall s>n/2.
  $$
  In particular: 
  \begin{align*}
&|u^\e|^2 \Tend \e 0 |a|^2\quad \text{in } L^\infty_TH^s, \text{ and}\\
&\e\IM\(\overline{u^\e
     }  \nabla u^\e   \) \Tend \e 0
     |a|^2\nabla\(\eik+\phi\) \quad \text{in }L^\infty_TX^s,\ \forall s>n/2.
  \end{align*}
\end{theorem}
Recall that in general, none of the terms $a$ or $a^\e$ is in
$L^2(\R^n)$. Though, the difference $a^\e-a$ is in $L^2(\R^n)$, and
asymptotically small as $\e \to 0$.
Note that $(\rho,\bv):= (|a|^2,\nabla
(\phi+\eik))$ solves the Euler--Poisson system:
\begin{equation}
  \label{eq:eulerpoisson}
\left\{
  \begin{aligned}
  \d_t \rho +\nabla\cdot (\rho \bv)&=0,\\
\d_t \bv +\bv\cdot \nabla \bv +\nabla\ext + \nabla\poi&=0,\\
\Delta \poi = q\(\rho -c\)\quad &; \quad \nabla\poi(t,x)\Tend
  {|x|}{\infty}0,\ \poi(t,0)=0. 
\end{aligned}
\right.
\end{equation}
The existence of solutions to
\eqref{eq:eulerpoisson} under Assumption~\ref{hyp:lin} is new.
\smallbreak

This paper borrows several ideas from \cite{CaBKW}, 
\cite{PG05} and \cite{Grenier98}. As we have already mentioned, an
important difference with \cite{Grenier98} is that the underlying wave
equation associated to \eqref{eq:systcomplet} is semi-linear, and not
quasi-linear. The reduction to \eqref{eq:systcomplet2} is similar to
the approach in \cite{CaBKW}. Several important differences should be
pointed out. First, we work in Zhidkov spaces instead of Sobolev
spaces, an aspect which requires some extra care. Integrating the
Poisson equation, especially when we have $\Delta \poi^\e\in L^2(\R^n)$
and not necessarily $\Delta \poi^\e\in L^1(\R^n)$, is also a new
problem. Finally, the propagation of the initial assumption
$|a_0^\e|^2 -c \in L^2(\R^n)$ turns out to be different from the
phenomenon studied in \cite{PG05}. As we shall see in
Section~\ref{sec:eik}, the presence of quadratic ``geometric'' quantities
(such as an external harmonic potential) requires a highly non-trivial
adaptation of the approach in \cite{CaBKW}. 
\smallbreak

The rest of this paper is organized as follows. In
Section~\ref{sec:general}, we collect various technical estimates, 
in order not to interrupt the proofs later on. In
Section~\ref{sec:lin}, we prove
Theorem~\ref{theo:existence}. Theorem~\ref{theo:convergence} is proved
in Section~\ref{sec:convergence}. In Part~\ref{part:2}
(Sections~\ref{sec:eik}--\ref{sec:reduction}), we consider the case
when $c-1\in L^1\cap H^\infty$, and the external potential and the
initial phase 
contain quadratic terms. 
  In Part~\ref{part:3}
(Section~\ref{sec:integrable}), we 
assume $c\in L^1(\R^n)$, and prove a refined
convergence result. 

\begin{remark}
  Before leaving this introduction, let us explain why
we concentrated on the whole space problem. Indeed, some problems
require considering the periodic case (see \cite{SP-Per} and the
references therein), where the space variable belongs to the torus
$\T^n$. As a matter of fact, the periodic case is easier. 
This follows from two observations: first, the 
computations below
apply {\em mutatis mutandis} in the periodic setting; and second, for
all $\sigma\in\R$,
the operator $\Delta^{-1}\nabla$ is well-defined in
$H^\sigma(\T^n)$. 
\end{remark}


\section{Estimates in Lebesgue, Sobolev and Zhidkov spaces}\label{sec:general}

This section serves as the requested background for what follows. The
proofs of easy or classical results are left out.
We first recall a consequence of the Hardy-Littlewood-Sobolev inequality,
which can be found in \cite[Th.~4.5.9]{Hormander1} or \cite[Lemma~7]{PG05}:
\begin{lemma}\label{lem:hormander} 
If $\varphi\in {\mathcal D}'(\R^n)$ is such that $\nabla \varphi\in
L^p(\R^n)$ for some $p\in ]1,n[$, then there 
exists a constant $\gamma$ such that $\varphi-\gamma \in L^q(\R^n)$, with
$1/p=1/q+1/n$. 
\end{lemma}
This shows that under Assumption~\ref{hyp:lin}, the doping
profile is of the form
\begin{equation*}
  c = \gamma + \widetilde c,\quad\text{where }\gamma\text{ is a
  constant, and }\widetilde c \in L^{\frac{2n}{n-2}}(\R^n), \nabla
  \widetilde c  \in H^\infty.
\end{equation*}
Define the Fourier transform as
\begin{equation*}
  \F \varphi(\xi)= \widehat \varphi(\xi
  )=\frac{1}{(2\pi)^{n/2}}\int_{\R^n}e^{-ix\cdot\xi} \varphi(x)dx. 
\end{equation*}
\begin{lemma}\label{lem:sobolev}
  Let $n\ge 3$. For every $s>n/2$, there exists 
  $C_{s}$ such that  
  \begin{equation}\label{esti:sobolev}
    \|\varphi\|_{L^\infty}\le C_{s}\|\nabla \varphi\|_{H^{s-1}},
    \quad 
    \forall \varphi \in H^{s}(\R^n).
  \end{equation}
\end{lemma}
\begin{remark*}
In space dimension $n\le 2$, low frequencies rule out the above
inequalities. For instance, in space dimension $n=1$, the function 
\begin{equation*}
  f(x) = \int_0^x\frac{dy}{\sqrt{1+y^2}}=\operatorname{arg}\sinh (x)
\end{equation*}
is not in $L^\infty(\R)$, but its derivative is in $H^\infty$.
In space dimension $n=2$, consider the
function 
\begin{equation*}
  f(x_1,x_2) =\log \left|\log (x_1^2 +x_2^2)\right|.
\end{equation*}
One can check that $\nabla f\in H^\infty$, while clearly,
$f\not\in L^\infty(\R^2)$. 
\end{remark*}
\noindent {\bf Warning} (Homogeneous Sobolev spaces). 
It may be
tempting to restate Lemma~\ref{lem:sobolev} in terms of homogeneous
Sobolev spaces. Recall that, for $s>0$, the homogeneous Sobolev space
$\dot H^s$ is 
defined as the completion of the Schwartz space ${\mathcal S}(\R^n)$
for the norm 
\begin{equation*}
  \|\varphi\|_{\dot H^s}= \||\xi|^s \widehat \varphi\|_{L^2}. 
\end{equation*}
More precisely, one might want to replace the right-hand side of
\eqref{esti:sobolev} 
with $\|\varphi\|_{\dot H^1}+ \|\varphi\|_{\dot H^s}$ and consider
$\varphi \in \dot H^1 \cap \dot H^s$ only. This is extremely delicate,
since $\dot H^s$ is not a Hilbert space when $s\ge n/2$. 

\begin{lemma}\label{lem:sobolevp}
  Let $n\ge 3$, $q\ge 2$ and $s>n/2-1$. There exists $C=C(n,q,s)$
  such that for all $\varphi\in L^q(\R^n)$ with $\nabla \varphi
  \in H^s(\R^n)$, 
  \begin{equation*}
    \|\varphi\|_{L^\infty} \le C\( \|\varphi\|_{L^q} + \|\nabla
    \varphi\|_{H^s}\). 
  \end{equation*}
\end{lemma}
\begin{proof}
  The usual Sobolev embedding yields, for any $\sigma >n/q$,
  \begin{equation*}
   \|\varphi\|_{L^\infty} \les \|\varphi\|_{L^q} + \||\nabla|^\sigma
    \varphi\|_{L^q}.  
  \end{equation*}
On the other hand, for $k=n(1/2-1/q)$,
\begin{equation*}
    \||\nabla|^\sigma \varphi\|_{L^q}\les \||\nabla|^\sigma
    \varphi\|_{H^k} \les \|\nabla
    \varphi\|_{H^{k+\sigma -1}}, 
  \end{equation*} 
provided that $\sigma \ge 1$.
If $s>n/2-1$,  $\sigma$ given by $s=n(1/2-1/q)+\sigma
-1$ is such that  $\sigma >n/q$ and $\sigma \ge 1$. The above two
estimates then yield the lemma. 
\end{proof}
\begin{lemma}\label{lem:poisson}
  Let $n\ge 3$. For every $s\ge 0$, 
  $\nabla \Delta^{-1}$ maps $L^1(\R^n)\cap H^s(\R^n)$ to
  $H^{s+1}(\R^n)$:  there exists
  $C_s$ such that 
  \begin{equation*}
    \|\nabla \Delta^{-1} \varphi\|_{H^{s+1}}\le C_s \(\|\varphi\|_{L^1}
    + \|\varphi\|_{H^s}\), \quad \forall \varphi \in L^1(\R^n)\cap
    H^s(\R^n). 
  \end{equation*}
\end{lemma}
The following variant of the classical Kato-Ponce estimates can be found in 
\cite[Theorem~5]{LannesJFA}:
\begin{lemma}
Let $n\ge 1$ and $s>n/2+1$. Denote $\Lambda = (\id-\Delta)^{1/2}$. There
exists a constant $C_{s}$ such that, 
for all $f\in X^{s+1}(\R^n)$ and all $u\in H^{s-1}(\R^{n})$, 
\begin{equation}\label{Kato-Ponce}
\left\lVert f \Lambda^{s} u - \Lambda^s(fu) \right\rVert_{L^2} \le
C_{s}\(\left\lVert \nabla f\right\rVert_{L^\infty}\left\lVert
  u\right\rVert_{H^{s-1}}+ \left\lVert \nabla
  f\right\rVert_{H^{s-1}}\left\lVert
  u\right\rVert_{L^\infty}\). 
\end{equation}
\end{lemma}

\begin{lemma}\label{lem:product}
Let $s>n/2$. The Sobolev space $H^s(\R^n)$ and the Zhidkov space
$X^s(\R^n)$ are algebras: there exists
a constant $C_{s}$ such that, for all $u,v\in
H^s(\R^n)$  and $a,b\in X^s(\R^n)$,
\begin{equation*}
\left\lVert u v \right\rVert_{H^s}\le C_{s}\left\lVert
  u\right\rVert_{H^s}\left\lVert v\right\rVert_{H^s}\quad ;\quad
  \left\lVert a b \right\rVert_{X^s}\le C_{s}\left\lVert 
  a\right\rVert_{X^s}\left\lVert b\right\rVert_{X^s}.
\end{equation*}
There exists $C_{s}$ such that for all  $v\in
H^s(\R^n)$ and $a\in 
X^s(\R^n)$,
\begin{equation*}
  \left\| av\right\|_{H^s} 
+\left\| av\right\|_{X^s} \le C_s \left\|
  v\right\|_{H^s}\left\| a\right\|_{X^s} .  
\end{equation*}
There exists $C_{s}$ such that for all  $a\in 
X^{s}(\R^n)$ and $b\in 
X^{s+1}(\R^n)$,
\begin{equation*}
  \left\lVert a \nabla b \right\rVert_{H^s}\le C_{s}\left\lVert 
  a\right\rVert_{X^s}\left\lVert b\right\rVert_{X^{s+1}}.
\end{equation*}
\end{lemma}

In order to use Arzela--Ascoli's theorem, we will invoke: 
\begin{lemma}\label{lem:compacite}
  Let $\sigma >n/2$ and $(\varphi_j)_{j\in\N}$ be a bounded  sequence
  in $X^\sigma(\R^n)$. For all $\sigma'<\sigma$, 
  there exists a subsequence which converges in $H^{\sigma'}_{\rm
  loc}(\R^n)$. 
\end{lemma}
\begin{proof}
This follows from the fact that, for all test
 function $\chi\in C_{0}^\infty(\R^n)$, $(\chi\varphi_j)_{j\in\N}$ is
 a bounded  sequence in $H^\sigma(\R^n)$. 
\end{proof}
\begin{remark*}
  It might seem more natural to state a precompactness result in 
\begin{equation*}
    X^{\sigma'}_{\rm loc}(\R^n)\defn \left\{ \varphi \in L^\infty_{\rm
    loc}(\R^n)\ ;\ \nabla \varphi \in H^{\sigma'-1}_{\rm
    loc}(\R^n)\right\}\cdot
  \end{equation*}
Actually, one can check that for $\sigma'>n/2$, $X^{\sigma'}_{\rm
  loc}(\R^n)= H^{\sigma'}_{\rm loc}(\R^n)$.
\end{remark*}

\begin{lemma}\label{lem:infini}
Let $n\ge 3$ and $s>n/2$. 
  \begin{itemize}
  \item For all $p >\frac{2n}{n-2}$, there exists $C=C(s,p,n)$ such that:
  \begin{equation*}
    \|\Delta^{-1}\nabla f\|_{L^p} \le C \|f\|_{H^s},\quad \forall f\in
    H^s. 
  \end{equation*}
  \item There exists $C=C(s,n)$ such that:
    \begin{equation*}
     \left\|{\mathcal F}\(\Delta^{-1}\nabla
     f\)\right\|_{L^1} \le C   \|f\|_{H^s}, \quad \forall f\in H^s.
    \end{equation*}
  \end{itemize}
\end{lemma}
\begin{proof}
Essentially, we use the property $\widehat f\in L^2$ for low
frequencies, and  $\widehat f\in L^1$ for high frequencies (${\mathcal
  F}(H^s)\subset L^1$ if $s>n/2$). For $ p>2n/(n-2)$, 
\begin{align*}
  \left\|
  \frac{\xi}{|\xi|^2}\widehat f\right\|_{L^{p'}}
\lesssim \left\| |\xi|^{-1}\right\|_{L^{\frac{2p'}{2-p'}}(|\xi|<1)}\left\|
  \widehat f\right\|_{L^{2}} + \left\|
  |\xi|^{-1}\right\|_{L^\infty(|\xi|\ge 1)}\left\|
  \widehat f\right\|_{L^{p'}}. 
\end{align*}
The norms involving $|\xi|^{-1}$ are finite since $p>2n/(n-2)$. For
$s>n/2$,  
\begin{equation*}
  \left\| \widehat f\right\|_{L^{1}}\lesssim \left\| \<\xi\>^{s}\widehat
  f\right\|_{L^{2}}= \|f\|_{H^s}.  
\end{equation*}
The first point follows from the Hausdorff--Young inequality:
\begin{equation*}
  \|\Delta^{-1}\nabla f\|_{L^p} \lesssim \left\|
  \frac{\xi}{|\xi|^2}\widehat f\right\|_{L^{p'}}
\end{equation*}
The second point is straightforward, with $p'=1$.
\end{proof}  


\part{Sublinear eikonal phase}\label{part:1}


\section{Proof of Theorem~\ref{theo:existence}}
\label{sec:lin}

Our first task is to construct a solution to
\eqref{eq:systcomplet2}. As explained in the introduction, it is
convenient to introduce 
the ``velocity'' $v^\e=\nabla \phi^\e$. Denoting $v_{\ei}=\nabla
\eik$, and recalling that $v_{\ei}$ is a function of time only, we infer 
from \eqref{eq:systcomplet2} that $(a^\e,v^\e)$ has to solve:
\begin{equation}\label{eq:systav0}
\left\{
\begin{aligned}
 & \d_t v^\e + \(v_{\ei}+v^\e\)\cdot
  \nabla v^\e+\nabla V_{\rm pert} +
  \nabla\poi^\e=0,\\
&\d_t a^\e +\(v_{\ei}+v^\e\) \cdot \nabla a^\e +\frac{1}{2}a^\e
  \nabla\cdot v^\e
  = i\frac{\e}{2}\Delta a^\e,\\ 
 &\Delta \poi^\e =q\big(|a^\e|^2-c\big),
\end{aligned}
\right.
\end{equation}
together with
\begin{equation}\label{ci:systav0}
\nabla \poi^\e(t,x)\Tend
  {|x|}{\infty}0\quad ;\quad \poi^\e(t,0)=0\quad ;\quad v^\e_{\mid
  t=0}=\nabla\phi_0 \ ;\ a^\e_{\mid t=0}=a_0^\e . 
\end{equation}
In the context of Assumption
\ref{hyp:lin}, we show  
that the solutions of \eqref{eq:systav0}--\eqref{ci:systav0} exist and
are uniformly bounded for a time interval independent of $\e$. 
\begin{proposition}\label{prop:EPoisson0} 
Let Assumption~\ref{hyp:lin} be satisfied. Let $s>n/2$. 
For all $ M>\ M_{0}>0$, there exists $T>0$
such that, if for all $\e\in 
[0,1]$,
\begin{equation}\label{TM:CI0}
\left\lVert \nabla
  \phi_{0}\right\rVert_{H^{s+2}}+ 
\bigl\lVert \left\lvert a_{0}^{\e} \right\rvert
^2-c\bigr\rVert_{L^{2}}+\left\lVert 
  a_{0}^{\e}\right\rVert_{X^{s+1}}\le  M_{0}, 
\end{equation}
then the Cauchy problem~\eqref{eq:systav0}--\eqref{ci:systav0} has a
unique classical solution   
$(v^{\e},a^{\e})$ in $ C^{\infty}([0,T]\times\R^{n})$ such that
\begin{equation}\label{TH:NORM0}
\left\lVert  v^{\e}\right\rVert_{L^\infty_T X^{s+2}}
+ \bigl\lVert \left\lvert a^{\e} \right\rvert
^2-c\bigr\rVert_{L^\infty_T L^{2}} 
+ \left\lVert  a^{\e}\right\rVert_{L^\infty_TX^{s+1}}\le M.
\end{equation}
\end{proposition}

As suggested by the above statement, we
construct $\nabla \poi^\e$ (only the gradient of $\poi^\e$ is present
in \eqref{eq:systav0}), and the condition $\poi^\e(t,0)=0$ is given
only to insure uniqueness for $\poi^\e$ (even though it is not
stated in the above result). Therefore, we shall neglect this
condition for a while.

\subsection{Regularized equations}
%
Let $\jmath$ be a $C^\infty$ function of $\xi\in\R^n$, with
$$
0\le \jmath\le 1,\quad \jmath(\xi)=1 \text{  for  }  \left\lvert
  \xi\right\rvert  \le 1,  \quad  
\jmath(\xi)=0 \text{  for  }\left\lvert \xi\right\rvert  \ge 2, \quad
\jmath(\xi)=\jmath(-\xi). 
$$              
Set $\jmath_{\para}(\xi)\defn \jmath(\para\xi)$, for $\para >0$ and
$\xi\in\R^n$;   
$\jmath_\para$ is supported in the ball of radius $2/h$ about the
origin. Define $J_\para$ as the Fourier multiplier with symbol
$\jmath_\para$:  
$$
J_\para\defn \jmath(\para D_x).
$$
Also, for our purpose it is interesting to introduce a family 
of operators that cut the low frequency component of a
function. Indeed, the Poisson term  
$\nabla V_{\rm p}^\e =
q\Delta^{-1}\nabla(\left\lvert a^\e\right\rvert
^2-c)$,  
is not well defined in general. We replace the operator
$q\Delta^{-1}\nabla$  
by a family of operators $R_\para \nabla$ well defined on Sobolev
spaces and prove that, in the end, there is no need  
to estimate the low frequency component of $\nabla V_{\rm p}^\e$. 
To do that, we set 
$$
 G_\para=\id - J_{1/\para},
$$
that is, $G_\para$ is the Fourier multiplier with symbol
$1-\jmath_{1/\para}$, which is  
supported in $\{ \left\lvert \xi\right\rvert \ge
\para\}$. Consequently, the operator 
$$
R_\para\defn q \Delta^{-1} G_\para,
$$
is bounded in all Sobolev spaces (with operator norm going to
$+\infty$ when $h$ tends to $0$). More precisely, there exists a constant $C$
such that, for all  $\sigma\ge 0$, 
$$
\left\lVert \Delta^{-1}
  G_\para\right\rVert_{H^\sigma\rightarrow H^{\sigma+2}}\le 
C \para^{-2}. 
$$
Consider the following approximation
of \eqref{eq:systav0}:
\begin{equation}\label{eq:systav02}
\left\{
\begin{aligned}
  &\d_t v_\para^\e + J_\para\(\(v_{\ei}+v_\para^\e\)\cdot\nabla J_\para
  v_\para^\e\)  
  + \nabla V_{\rm pert} =
  -R_\para \nabla(\left\lvert a_\para^\e\right\rvert ^2-c) ,\\
 &\d_t a_\para^\e+J_\para \(\(v_{\ei}+v_\para^\e\) \cdot \nabla J_\para
  a_\para^\e\) +\frac{1}{2} a_\para^\e 
 \nabla\cdot v_\para^\e =
  i\frac{\e}{2}\Delta J_\para^2 a_\para^\e. 
\end{aligned}
\right.
\end{equation}
We keep the same initial data:
\begin{equation}\label{ci:systav02}
{v_{\para}^{\e}}_{\arrowvert t=0}=\nabla\phi_0 \ ;\
{a_{\para}^{\e}}_{\mid t=0}= a_0^\e. 
\end{equation}
Note that Assumption~\ref{hyp:lin} implies that
${v_{\para}^{\e}}_{\arrowvert t=0}$ is in $H^\infty$ and is
independent of $\e\in [0,1]$ and $\para\in ]0,1]$, while
${a_{\para}^{\e}}_{\mid t=0}$ is in $X^\infty$, and uniformly bounded
in $X^s$ for any $s>n/2$, for $\e\in [0,1]$ and $\para\in ]0,1]$.

The point is that the regularized equations
\eqref{eq:systav02}--\eqref{ci:systav02} 
have been chosen so that  the Cauchy problem can be solved as in the
standard framework of Sobolev spaces:
\begin{lemma}\label{lem:existEDO0}
Let $s>n/2$. For all $\e\in [0,1]$ and all $\para\in ]0,1]$ 
there exists $T_\para^\e>0$ such that the Cauchy problem
\eqref{eq:systav02}--\eqref{ci:systav02}  has a unique solution  
$(v_\para^\e,a_\para^\e)\in C^{1}([0,T_\para^\e];H^{s+2}(\R^n)\times
X^{s+1}(\R^n))$. 
\end{lemma}
\begin{proof}
The proof is based on the usual theorem for ordinary differential equations. 
Set 
$\bu_\para^\e=(v_\para^\e,a_\para^\e)$ and  
we rewrite \eqref{eq:systav02} under the form 
$$
\d_{t}\bu_\para^\e=F_1(\e,\para,\bu_\para^\e)+F_2(t)\bu_\para^\e+F_3(t,x),
$$
where $F_1(\e,\para,\bu)$ is at most quadratic in $\bu$, and we have
used the property that $ v_{\ei}$ is a function of 
time only.  We 
have to verify that the functions $F$ are smooth.  
This follows from Lemmas~\ref{lem:sobolev} and \ref{lem:product},  and
the fact that the operators   
$R_\para$ and $\Delta J_\para $ are of order $-2$ and
$0$ respectively:
\begin{align*}
&\left\lVert J_\para(v_\para^\e\cdot\nabla J_\para
  v_\para^\e)\right\rVert_{H^{s+2}} 
\les \left\lVert v_\para^\e\right\rVert_{H^{s+2}}\left\lVert \nabla
  J_\para v_\para^\e\right\rVert_{H^{s+2}} 
\les \para^{-1}\left\lVert v_\para^\e\right\rVert_{H^{s+2}}^2,\\
&\left\lVert J_\para(v_\para^\e \cdot \nabla J_\para
  a_\para^\e)\right\rVert_{X^{s+1}} \les  
 \left\lVert v_\para^\e\right\rVert_{H^{s+1}}\left\lVert \nabla
  J_\para a_\para^\e\right\rVert_{X^{s+1}} 
\les \para^{-1}\left\lVert v_\para^\e\right\rVert_{H^{s+1}}\left\lVert
  a_\para^\e\right\rVert_{X^{s+1}},\\ 
&\left\lVert a_\para^\e\nabla \cdot v_\para^\e\right\rVert_{X^{s+1}}\les
\left\lVert v_\para^\e\right\rVert_{H^{s+2}}\left\lVert
  a_\para^\e\right\rVert_{X^{s+1}},\\ 
&\bigl\lVert R_\para \nabla \left\lvert a_\para^\e\right\rvert
  ^2\bigr\rVert_{H^{s+2}} 
\les\para^{-2} \bigl\lVert \nabla \left\lvert a_\para^\e\right\rvert
  ^2\bigr\rVert_{H^{s}}  
\les \para^{-2} \bigl\lVert  a_\para^\e\bigr\rVert_{X^{s+1}}^2, \\
&\left\lVert \Delta J_\para^2 a_\para^\e\right\rVert_{X^{s+1}}
= \left\lVert J_\para\Delta J_\para a_\para^\e\right\rVert_{X^{s+1}}
\le\left\lVert \Delta J_\para a_\para^\e\right\rVert_{H^{s+1}}
\les\para^{-2}\left\lVert a_\para^\e\right\rVert_{X^{s+1}} .
\end{align*}
\end{proof}

\subsection{Uniform bounds}
\label{sec:uniflin}

To prove Proposition~\ref{prop:EPoisson0}, the analysis of 
\eqref{eq:systav02} contains at least two parts:  
first, an existence and uniform boundedness result for a time
independent of the small parameters $\e$ and $\para$; and  
second, a convergence result when  $\para \to 0$.
Here, we prove that the solutions $(v_\para^\e,a_\para^\e)$ exist and
they are uniformly bounded for  
a time independent of the parameters $\e$ and $\para$. 
Below, $T_\para^{\e *}$ denotes the lifespan, that is the supremum of
all the positive times $T_\para^\e$ such that the  
Cauchy problem for \eqref{eq:systav02}--\eqref{ci:systav02} has a
unique solution in  $C^{1}([0,T_\para^\e];H^{s+2}(\R^n)\times
X^{s+1}(\R^n))$. 
\begin{proposition}\label{prop:gronwall0}
Let $s>n/2$. There exists a continuous function $g\colon
\R_+^*\rightarrow \R_+^*$ such that, for all $\e\in [0,1]$ and all  
$\para\in ]0,1]$, the norm $M_\para^\e\colon
[0,T_\para^{\e*}]\rightarrow \R_+^*$ defined by 
$$
M_\para^\e(T)\defn \| a_\para^\e\|_{L^\infty_TX^{s+1}} + 
\left\| |a_\para^\e|^2-c\right\|_{L^\infty_T L^2}
+ \|\nabla v_\para^\e\|_{L^\infty_T H^{s+1}}, 
$$
satisfies the estimate: 
$M_\para^\e(T)\le M_\para^\e(0)e^{Tg( M_\para^\e(T))}$, $\forall T\in
[0,T_\para^{\e *}]$. 
\end{proposition}
\begin{proof}
Before we proceed, two comments are in order. Firstly, 
the functions $(v_\para^\e,a_\para^\e)$ are smooth ($C^1$ in time with
values in Sobolev/Zhidkov spaces), so that  
it is easily verified that all the following computations are meaningful. 
Secondly, it is useful to note that, in view of
Lemma~\ref{lem:sobolev}, it suffices to prove that  
\begin{equation}\label{ubound0}
m_\para^\e(T)\le m_\para^\e(0)e^{Tg( \bound(T))}, \qquad \forall T\in
[0,T_\para^{\e *}], 
\end{equation}
where
\begin{equation*}
  m_\para^\e(t) = \| \nabla a_\para^\e\|_{L^\infty_TH^{s}} + 
\left\| |a_\para^\e|^2-c\right\|_{L^\infty_T L^2}
+ \|\nabla v_\para^\e\|_{L^\infty_T H^{s+1}},
\end{equation*}
and
$$
\bound(T)\defn M_\para^\e(T) + \left\lVert
  v_\para^\e\right\rVert_{L^{\infty}([0,T]\times \R^n)}. 
$$
Indeed, Lemma~\ref{lem:sobolev} provides us with a constant  $C_s$ such that 
$\bound\le C_s M_\para^\e $, and we have:
\begin{lemma}\label{lem:contlinfini0}
  Let $s>n/2$ and $c\in X^\infty$. There exists a constant $K$ such
  that, for all 
  $T>0$ and $\varphi\in X^\infty$,
  \begin{equation*}
    \|\varphi\|_{L^\infty(\R^n)}^2 \le K 
    \left\lVert\lvert\varphi\rvert^2-c\right\rVert_{L^2}+K
    \|\nabla \varphi
    \|_{H^{s}}^2+K \|c\|_{X^{s+1}}. 
  \end{equation*}
\end{lemma}
\begin{proof}
By Lemma~\ref{lem:sobolevp}, we have:
\begin{equation*}
  \left\||\varphi|^2-c\right\|_{L^\infty}
  \les \left\||\varphi|^2-c\right\|_{L^2 } +
  \left\|\nabla\bigl(|\varphi|^2-c\bigr)\right\|_{H^{s}}.  
\end{equation*}
Since $s>n/2$, 
\begin{align*}
  \left\lVert\nabla\bigl(\lvert\varphi\rvert^2-c\bigr)\right\|_{H^{s}}
  \les \left\|\varphi\right\|_{L^\infty}
  \left\|\nabla\varphi\right\|_{H^{s}}
  +\left\|\nabla\varphi\right\|_{H^{s}}^2 
  +  \left\|\nabla c\right\|_{H^{s}}. 
\end{align*}
Triangle inequality yields
\begin{align*}
  \left\|\varphi\right\|_{L^\infty}^2 \les &\ 
  \left\| |\varphi|^2-c\right\|_{L^2 } +
  \left\|\varphi\right\|_{L^\infty} 
  \left\|\nabla\varphi\right\|_{H^{s}}
  +\left\|\nabla\varphi\right\|_{H^{s}}^2 
  +  \left\|c\right\|_{X^{s+1}},
\end{align*}
hence, the desired result follows by Young's inequality.
\end{proof}

With these preliminaries established, to prove \eqref{ubound0}, we
begin by estimating the $L^2$ norm of 
$\lvert a_\para^\e \rvert ^2-c$. To do that, we start from
\begin{equation*}
  \frac{d}{dt}\bigl\lVert\left\lvert
  a_{\para}^\e\right\rvert ^2-c\bigr\rVert_{L^2}^2 
  \le 2 \bigl\lVert\lvert 
  a_{\para}^\e\rvert ^2-c\bigr\rVert_{L^2}\bigl\lVert\d_t 
  \bigl(\lvert  
  a_{\para}^\e\rvert ^2-c\bigr)\bigr\rVert_{L^2}. 
\end{equation*}
The second factor in the right hand side is estimated by
\begin{equation*}
\bigl\lVert\d_{t}\bigl(\left\lvert a_{\para}^\e\right\rvert
  ^2-c\bigr)\bigl\rVert_{L^2}
\le 2\left\lVert a_{\para}^\e\right\rVert_{L^{\infty}}
\left\lVert
  \d_{t}a_{\para}^\e\right\rVert_{L^2}. 
\end{equation*}
Directly from the equations, we find that for bounded times,
\begin{equation*}
\left\lVert \d_t a_\para^\e \right\rVert_{L^2}\les 
\(1+\left\lVert v_\para^\e
   \right\rVert_{L^\infty}\)\left\lVert\nabla a_\para^\e
   \right\rVert_{L^2} +\left\lVert a_\para^\e \right\rVert_{L^\infty}
   \left\lVert\nabla v_\para^\e \right\rVert_{L^2}
+ \left\lVert \Delta a_\para^\e  \right\rVert_{L^2}. 
\end{equation*}
Consequently, we obtain
\begin{equation}
  \label{eq:borneL20}
  \frac{d}{dt}\bigl\lVert\lvert
  a_{\para}^\e\rvert ^2-c\bigr\rVert_{L^2}^2 \les M^\e_\para 
  \( 1+M^\e_\para \)^2.
\end{equation}
We now turn to the estimate of the $H^s$ norm of $\nabla a_\para^\e$. 
Set $Q\defn \Lambda^{s}\nabla$, where $\Lambda =(\id-\Delta)^{1/2}$. 
Since $[\nabla,Q]=0=[J_\para,Q]$, by commuting $Q$ with the equation for
$a_\para^\e$, we find: 
\begin{equation*}
\d_t Q a_\para^\e 
+J_\para\(\(v_{\ei}+v_\para^\e\) \cdot \nabla J_\para Q a_\para^\e\)
  -i\frac{\e}{2}\Delta J_\para^2Qa_\para^\e =f_\para^\e,
\end{equation*}
with
\begin{align*}
f_\para^\e\defn J_\para \([v_\para^\e,Q]\cdot\nabla J_\para
a_\para^\e\)-\frac{1}{2} Q\bigl(a_\para^\e\nabla\cdot v_{\para}^{\e}\bigr).
\end{align*}
Notice that $J_\para$ is self-adjoint. We use the
following convention for the scalar product in $L^2$:
\begin{equation*}
  \scal{f}{g}\defn \int_{\R^n} f(x)\overline{g(x)}dx.
\end{equation*}
We have, since $\nabla
v_{\ei}\equiv 0$: 
\begin{align*}
\RE\scal{i\Delta J_\para^2
  Qa_\para^\e}{Qa_\para^\e}&=\RE\scal{i\Delta J_\para
  Qa_\para^\e}{J_\para Qa_\para^\e}=0,\\ 
2\scal{J_\para((v_{\ei}+v_\para^\e)\cdot\nabla J_\para Qa_\para^\e)}{Q
  a_\para^\e}&= 
2\scal{(v_{\ei}+v_\para^\e)\cdot\nabla J_\para Qa_\para^\e}{J_\para
  Qa_\para^\e}\\ 
&=-\scal{(\nabla\cdot v_\para^\e)J_\para Qa_\para^\e}{J_\para Qa_\para^\e}.
\end{align*}
Therefore,
\begin{align*}
\frac{d}{dt}\left\lVert
  Qa_\para^\e\right\rVert_{L^2}^2&=2\RE\scal{\d_t Q
  a_\para^\e}{Qa_\para^\e}\\ 
&= \RE\scal{(\nabla \cdot v_\para^\e)J_\para Q a_\para^\e}{J_\para
  Qa_\para^\e}+2\RE\scal{f_\para^\e}{Qa_\para^\e}\\ 
&\le \left\lVert \nabla v_\para^\e \right\rVert_{L^\infty}\left\lVert
  J_\para Q a_\para^\e\right\rVert_{L^2}^2+2\left\lVert
  f_\para^\e\right\rVert_{L^2}\left\lVert
  Qa_\para^\e\right\rVert_{L^2}.
\end{align*}
We now have to estimate the $L^2$ norm of $f_\para^\e$. The first term
is estimated  
by way of the commutator estimate \eqref{Kato-Ponce} and the Sobolev
embedding: 
\begin{align*}
\left\lVert J_\para \([v_\para^\e,Q]\cdot\nabla J_\para a_\para^\e\)
\right\rVert_{L^2} 
&\les \left\lVert [v_\para^\e,Q]\cdot\nabla J_\para a_\para^\e
\right\rVert_{L^2} \\ 
&\les \(\left\lVert \nabla v_\para^\e \right\rVert_{L^\infty}
+\left\lVert \nabla v_\para^\e \right\rVert_{H^{s+1}} \)\left\lVert
  \nabla J_\para a_\para^\e\right\rVert_{H^s} \\ 
&\les  \left\lVert \nabla v_\para^\e \right\rVert_{H^{s+1}}\left\lVert
  \nabla a_\para^\e\right\rVert_{H^s}. 
\end{align*}
To estimate the last term, we use Lemma~\ref{lem:product}, to obtain
\begin{equation*}
\left\lVert Q \bigl(a_\para^\e\nabla\cdot
  v_\para^\e\bigr)\right\rVert_{L^2}
  \les  \left\lVert a_\para^\e\nabla\cdot
  v_\para^\e\right\rVert_{H^{s+1}}
\les  \left\lVert a_\para^\e \right\rVert_{X^{s+1}}\left\lVert \nabla
  v_\para^\e\right\rVert_{H^{s+1}}. 
\end{equation*}
We infer that 
 \begin{equation*}
 \left\lVert f_\para^\e \right\rVert_{L^2}\les \left\lVert \nabla
   v_\para^\e\right\rVert_{H^{s+1}} 
\left\lVert a_\para^\e\right\rVert_{X^{s+1}}. 
 \end{equation*} 
Therefore, we end up with
\begin{equation}
\label{eq:aHs0}
\frac{d}{dt}\|\nabla a_\para^\e  \|_{H^s}^2 \lesssim \| \nabla
v_\para^\e \|_{H^{s+1}} \|  a_\para^\e  \|_{X^{s+1}}^2. 
\end{equation}
The technique for estimating $\nabla v_\para^\e$ in $H^{s+1}$ is
similar. Indeed, the analysis establishing the previous estimate also
yields 
\begin{align*}
  \frac{d}{dt}\|\nabla v_\para^\e  \|_{H^{s+1}}^2& \lesssim \(1 +\|
  \nabla   v_\para^\e \|_{H^{s+1}} \)\| \nabla v_\para^\e
   \|_{H^{s+1}}^2  \\ 
&\quad +  \| \nabla^2 V_{\rm pert}\|_{H^{s+1}}^2 
+\| \nabla R_\para \nabla\bigl( \left\lvert
  a_\para^\e \right\rvert ^2-c \bigr)\|_{H^{s+1}}^2.
\end{align*}
Since $\nabla R_\para\nabla$ is uniformly bounded from 
$H^{s+1}$ to itself, we obtain
\begin{equation*}
 \frac{d}{dt}\|\nabla v_\para^\e  \|_{H^{s+1}}^2\lesssim 
\(1 +\| \nabla
  v_\para^\e \|_{H^{s+1}} \)\| \nabla v_\para^\e  \|_{H^{s+1}}^2 +  1
  +\left\lVert |a_\para^\e |^2-c \right\rVert_{H^{s+1}}^2. 
\end{equation*} 
Next, noting that
\begin{align*}
  \left\| |a_\para^\e |^2-c \right\|_{H^{s+1}} &\lesssim \left\|
  |a_\para^\e |^2-c \right\|_{L^2} 
+\|\nabla c \|_{H^s} + \|\nabla |a_\para^\e|^2\|_{H^s}\\
 & \lesssim \left\|
  |a_\para^\e |^2-c \right\|_{L^2} 
+1 + \|a_\para^\e\|_{L^\infty}^2 +\|\nabla a_\para^\e\|_{H^s}^2 \\
&\lesssim  \left\|
  |a_\para^\e |^2-c \right\|_{L^2} +1 +\|a_\para^\e\|_{X^{s+1}}^2, 
\end{align*}
we conclude that
\begin{equation}
  \label{eq:vHs0}
  \frac{d}{dt}\|\nabla v_\para^\e  \|_{H^{s+1}}^2\lesssim C\bigl(
  \bound \bigr) \|\nabla v_\para^\e  \|_{H^{s+1}}^2 +C(M_\para^\e ). 
\end{equation}
Summing over \eqref{eq:borneL20}, \eqref{eq:aHs0} and \eqref{eq:vHs0},
Gronwall lemma yields the uniform estimate \eqref{ubound0}.
\end{proof}
Lemma~\ref{lem:existEDO0} and  Proposition~\ref{prop:gronwall0} 
yield the following result: 
\begin{corollary}\label{cor:borneunif0}
Let Assumption~\ref{hyp:lin} be satisfied, and let $s>n/2$. For
all $M>M_0>0$,  
there exists $T>0$ such that, if for all $\e\in
[0,1]$,  
\begin{equation*}
\left\lVert \nabla \phi_{0}\right\rVert_{H^{s+2}(\R^{n})}+
\bigl\lVert \left\lvert a_{0}^{\e} \right\rvert
^{2}-c(\cdot)\bigr\rVert_{L^{2}(\R^{n})}+\left\lVert 
  a_{0}^{\e}\right\rVert_{X^{s+1}(\R^{n})}\le M_{0}, 
\end{equation*}
then the Cauchy problem \eqref{eq:systav02}--\eqref{ci:systav02} has a
unique classical solution $(v_\para^\e,a_\para^\e)\in
C^{1}([0,T];H^{s+2}(\R^n)\times 
X^{s+1}(\R^n))$ satisfying  
\begin{equation*}
\left\lVert \nabla v_\para^{\e}\right\rVert_{L^\infty_T H^{s+1}}
+ \bigl\lVert  \left\lvert a_\para^{\e}
\right\rvert^{2}-c\bigr\rVert_{L^\infty_T L^{2}}  
+ \left\lVert  a_\para^{\e}\right\rVert_{L^\infty_T X^{s+1}}\le M.
\end{equation*}
\end{corollary}
\begin{remark}
  Refining the above computations thanks to Moser's calculus and tame
  estimates, we 
  can see that the above existence time $T$ can be taken independent
  of $s>n/2$ (see e.g. \cite[Section~2.2]{Majda} or
  \cite[Section~16.1]{Taylor3}). This explains why we did not
  emphasize its dependence 
  upon $s$, and why we consider different values for $s$
  below, without changing the notation $T$. 
\end{remark}

\subsection{Convergence of the scheme}
\label{sec:cvschemelin}
We first claim that $\d_t  v_\para^{\e}$ and
$\d_t a_\para^\e$ are bounded in $C([0,T];X^{s-1})$, uniformly for
$\para \in]0,1]$. To see this, 
by using Lemma~\ref{lem:sobolev}, \eqref{eq:systav02}  and
Corollary~\ref{cor:borneunif0}, the point is to verify that the term
$R_{\para}\nabla \( |a^\e_\para|^2 -c\)$, in the equation  
for $\partial_{t}v^\e_{\para}$, is uniformly bounded 
in $C([0,T];L^{\infty})$. Denote 
\begin{equation*}
W_\para^\e \defn R_\para \nabla\( |a^\e_\para|^2 -c\).
\end{equation*}
From Corollary~\ref{cor:borneunif0}, $W_\para^\e \in
C([0,T];H^{s+2})$, and $\nabla W_\para^\e $ is bounded in
$C([0,T];H^{s+1})$. In particular, Lemma~\ref{lem:sobolev} shows that
$W_\para^\e $ is bounded in $C([0,T];L^\infty)$.  
\smallbreak

From Lemma~\ref{lem:compacite} and Arzela--Ascoli's Theorem, for a
subsequence $\para'$ of 
$\para$, 
\begin{equation}\label{eq:cvdistri}
  v^\e_{\para'}\to v^\e \text{ and } a^\e_{\para'}
  \to a^\e 
  \text{ in } C([0,T];H^{s'}_{\rm loc}),\text{ as }\para'\to 0,  
\end{equation}
for any $s'<s-1$. Moreover, we have $v^\e,a^\e \in C_{w}([0,T];
X^{s})$. We can then pass to the limit in all the terms in
\eqref{eq:systav02}, except possibly the Poisson term, that is, the
right hand side in the equation for $v_\para^\e$. 
\smallbreak

To claim that $(v^\e,a^\e)$ solves
\eqref{eq:systav0}--\eqref{ci:systav0}, we introduce the Poisson
potential
\begin{equation*}
  V^\e_\para \defn  q \Delta^{-1}G_\para \( |a^\e_\para|^2 -c\). 
\end{equation*}
Then \eqref{eq:systav02} can be rewritten as:
\begin{equation}\label{eq:02eq}
\left\{
\begin{aligned}
  &\d_t v_\para^\e + J_\para\(\(v_{\ei}+v_\para^\e\)\cdot\nabla J_\para
  v_\para^\e\) 
  + \nabla V_{\rm pert} +\nabla V^\e_\para=0,\\
 &\d_t a_\para^\e+J_\para \(\(v_{\ei}+v_\para^\e\) \cdot \nabla J_\para
  a_\para^\e\) +\frac{1}{2} a_\para^\e 
 \nabla\cdot v_\para^\e=
  i\frac{\e}{2}\Delta J_\para^2 a_\para^\e,\\
&\Delta  V^\e_\para = q G_\para \( |a^\e_\para|^2 -c\).
\end{aligned}
\right.
\end{equation}
A  subsequence of $W_\para^\e$ converges in ${\mathcal D}'$ to some $W^\e\in
L^\infty([0,T]\times\R^n)$. Since 
$\nabla \times W_\para^\e = 0$ for every $\para\in ]0,1]$, we
deduce that $\nabla \times W^\e = 0$. We infer that there exists
$\poi^\e$ such 
that $W^\e = \nabla \poi^\e$ (see e.g. \cite[Prop.~1.2.1]{JYClivre}), and
we note
\begin{equation}\label{eq:17h30}
\nabla^2 \poi^\e \in C_w([0,T]; H^s).
\end{equation}
On the other hand, Corollary~\ref{cor:borneunif0} and Fatou's lemma
imply that $|a^\e|^2-c \in L^\infty([0,T];L^2)$. To prove that
$(v^\e,a^\e)$ solves \eqref{eq:systav0}--\eqref{ci:systav0}, we now
just have to check that $\Delta \poi^\e -q \(|a^\e|^2-c\)=0$. We
proceed in two steps: first, we prove that this quantity is a function
of time only. Then, since it is in $ L^\infty([0,T];L^2)$, we conclude
that it is necessarily zero. We have
\begin{equation*}
  \left\| \nabla\(\Delta \poi^\e -q \(|a^\e|^2-c\)\)\right\|_{L^2} \le
  \liminf_{\para \to 0} \left\| \nabla\(\Delta V^\e_\para -q
  \(|a^\e_\para|^2-c\)\)\right\|_{L^2} .
\end{equation*}
The last quantity is equal to: 
\begin{equation*}
  |q| \left\| \nabla J_{1/\para}
  \(|a^\e_\para|^2-c\)\right\|_{L^2}.
\end{equation*}
This goes to zero with $\para$, since
$|a^\e_\para|^2-c$ is uniformly bounded in $L^\infty_TL^2$: 
\begin{align*}
  \left\| \nabla J_{1/\para}
  \(|a^\e_\para|^2-c\)\right\|_{L^2} \les \left\|\xi  \jmath\(
  \frac{\xi}{\para} \) {\mathcal F}\(|a^\e_\para|^2-c\)\right\|_{L^2}
\les \para \left\|{\mathcal F}\(|a^\e_\para|^2-c\)\right\|_{L^2}
  \les \para.
\end{align*}
We infer that 
\begin{equation*}
 \nabla\(\Delta \poi^\e -q \(|a^\e|^2-c\)\)\equiv 0, 
\end{equation*}
that is, $\Delta \poi^\e -q \(|a^\e|^2-c\)$ is a function of time
only. We conclude that $(v^\e,a^\e)$ solves
\eqref{eq:systav0}--\eqref{ci:systav0}. 
\smallbreak

We prove additional regularity for $(v^\e,a^\e)$ by
showing that $(v^\e_\para - v^\e, a^\e_\para -a^\e)$ (and not a
subsequence) goes to zero in $L^\infty([0,T];X^{s+2}\times
X^{s+1})$. We will use:
\begin{lemma}\label{lem:cvop}
  Let $\varphi \in H^1$. Then:
  \begin{enumerate}
  \item $\displaystyle \left\| J_{1/\para}\varphi\right\|_{L^2}\to 0$
  as $\para \to 0$. 
  \item $\displaystyle \left\| (\id-J_{\para})\varphi\right\|_{L^2} +
  \left\| (\id-J_{\para}^2)\varphi\right\|_{L^2} \le 2 h
  \|\nabla \varphi\|_{L^2}$.
  \item There exists $C>0$ such that for all $\para\in]0,1]$, $\|R_\para
  \nabla^2\|_{L^2\to L^2}\le C$.
  \end{enumerate}
\end{lemma}
\begin{remark}
   Note that in the first point, $\varphi$ is supposed to be
  independent of $\para$. Otherwise, the conclusion needs not be true,
  which is easily checked by considering $\varphi_\para(x)=
  \para^{n/2}U(\para x)$, where $U\in {\mathcal S}(\R^n)$: $\left\|
  J_{1/\para}\varphi_{\para}\right\|_{L^2} = \| \jmath {\mathcal
  F}U\|_{L^2}$, is independent of $\para$. In the second point,
  $\varphi$ may of course depend on $\para$.   
\end{remark}
\begin{proof}
  For the first point, we write
  \begin{equation*}
    \left\| J_{1/\para}\varphi\right\|_{L^2} = \left\|
    \jmath\(\frac{\xi}{h} \) \widehat\varphi\right\|_{L^2},
  \end{equation*}
and we conclude with the Dominated Convergence Theorem. Next, we have:
\begin{equation*}
  \left\| (I-J_{\para})\varphi\right\|_{L^2} =\left\|
    \(1- \jmath\(h \xi \)\) \widehat\varphi\right\|_{L^2}\le
    h \left\|
    \(1- \jmath\(h \xi \)\)|\xi| \widehat\varphi\right\|_{L^2},
\end{equation*}
since the function $1- \jmath\(h \xi \)$ is supported in
$\{|\xi|>1/h\}$. The second term in the second point is treated
similarly. 
The last point follows from the fact that the symbol of the Fourier
multiplier $\Delta^{-1}\nabla^2$ is  bounded.  
\end{proof}
Denote
$(w_\para^\e,d^\e_\para)\defn (v^\e_\para - v^\e, a^\e_\para -a^\e)$,
and for $s>n/2+1$, introduce 
\begin{equation*}
  \rho^\e_\para (t) \defn \left\| w^\e_\para(t)\right\|_{L^\infty}+
  \left\| \nabla w^\e_\para(t)\right\|_{H^s} +\left\| d^\e_\para
  (t)\right\|_{H^s}. 
\end{equation*}
By construction, $\rho^\e_\para(0)=0$. Write the equation for
$(w_\para^\e,\nabla w_\para^\e,d^\e_\para)$ as: 
\begin{equation}\label{eq:reste}
\left\{
  \begin{aligned}
\d_t  w_\para^\e +(v_{\ei}+&v_\para^\e)\cdot \nabla w_\para^\e +w_\para^\e
    \cdot \nabla v^\e =\\
&=  -R_\para \nabla(\id -J_{1/\para})\( |a_\para^\e|^2 -
    |a^\e|^2\) +S_\para^\e,\\
    \d_t \nabla w_\para^\e +(v_{\ei}+&v_\para^\e)\cdot \nabla^2
    w_\para^\e +w_\para^\e 
    \cdot \nabla^2 v^\e +\nabla w_\para^\e\cdot \nabla v_\para^\e \\
+ \nabla
    v^\e\cdot \nabla w_\para^\e 
&
=  -R_\para \nabla^2(\id -J_{1/\para})\( |a_\para^\e|^2 -
    |a^\e|^2\) +\nabla S_\para^\e,\\
\d_t d_\para^\e + w_\para^\e\cdot \nabla a_\para^\e &+
    (v_{\ei}+v^\e)\cdot \nabla d_\para^\e 
    +\frac{1}{2}\( d_\para^\e \nabla \cdot v_\para^\e +a^\e \nabla \cdot
    w_\para^\e\)\\
&=  i\frac{\e}{2}\Delta
    d_\para^\e +\Sigma_\para^\e, 
   \end{aligned}
\right.
\end{equation}
where the source terms are given by
\begin{align*}
S_\para^\e =&v_\para^\e\cdot \nabla v_\para^\e - J_\para\(
  v_\para^\e\cdot \nabla J_\para v_\para^\e\)
 + q \Delta^{-1}\nabla J_{1/\para}\(|a^\e|^2-c\),\\
\Sigma_\para^\e =& v_\para^\e\cdot \nabla a_\para^\e - J_\para \(
  v_\para^\e\cdot \nabla J_\para a_\para^\e \) - i\frac{\e}{2} \Delta (\id
  -J_\para^2) a_\para^\e. 
\end{align*}
The error term $r_\para^\e$ may seem to involve too many quantities
(too much 
regularity), compared to the classical approach explained for instance
in \cite{Majda,Taylor3}. The usual approach would consist in
estimating $(w_\para^\e,d^\e_\para)$ in $L^2$ only. We cannot get such
estimates because of the Poisson term in $S_\para^\e$:  we can
prove it goes to zero in $X^s$, but not in $L^2$.  
We proceed in two steps:
\begin{enumerate}
\item We show that we can apply Gronwall lemma for $r_\para^\e(t)$,
  with sources terms  $S_\para^\e$ and
  $\Sigma_\para^\e$. 
\item We show that these source terms go to zero with $\para$ in the
  norms involved at the first step. 
\end{enumerate}
To estimate the first term of $r_\para^\e$, integrate in time the
first equation in \eqref{eq:reste}, and use
Corollary~\ref{cor:borneunif0}: 
\begin{align*}
 \|w_\para^\e(t)\|_{L^\infty} & \les \int_0^t \|\nabla
 w_\para^\e(\tau)\|_{L^\infty}d\tau + \int_0^t \|
 w_\para^\e(\tau)\|_{L^\infty}d\tau\\
+ \int_0^t& \left\|R_\para \nabla(\id -J_{1/\para})\( |a_\para^\e|^2 -
    |a^\e|^2\)(\tau)\right\|_{L^\infty} d\tau  +\int_0^t\| S_\para^\e
 (\tau)\|_{L^\infty}d\tau .
\end{align*}
Estimate the third term of the right hand side thanks to
Lemma~\ref{lem:sobolev}, Corollary~\ref{cor:borneunif0} and the last
point of Lemma~\ref{lem:cvop}:
\begin{align*}
 & \left\|R_\para \nabla(\id -J_{1/\para})\( |a_\para^\e|^2 - 
    |a^\e|^2\)(\tau)\right\|_{L^\infty} \\
    &\qquad \les 
    \left\| \nabla R_\para \nabla(\id -J_{1/\para})\( |a_\para^\e|^2 - 
    |a^\e|^2\)(\tau)\right\|_{H^{s-1}} \\
     & \qquad \les \left\|(\id
  -J_{1/\para})\( |a_\para^\e|^2 - |a^\e|^2\)(\tau)\right\|_{H^s}\\
& \qquad \les \left\|\( |a_\para^\e|^2 - 
    |a^\e|^2\)(\tau)\right\|_{H^s}\\
&\qquad \les \|(a_\para^\e -a^\e)(\tau)\|_{H^s}. 
\end{align*}
Using Sobolev embedding for the term in $\nabla w_\para^\e$, we end up with:
\begin{equation}
  \label{eq:resteLinf}
  \|w_\para^\e(t)\|_{L^\infty}  \les \int_0^t 
  r_\para^\e(\tau)d\tau + \int_0^t\| S_\para^\e
 (\tau)\|_{L^\infty}d\tau .
\end{equation}
Now estimate the $H^s$ norm of $\nabla w_\para^\e$. 
From the second equation in \eqref{eq:reste}, 
\begin{align*}
 \frac{d}{dt}\|\Lambda^s\nabla w_\para^\e\|_{L^2}^2 &= 2\RE
  \<\Lambda^s \nabla 
  w_\para^\e, 
  \d_t  \Lambda^s \nabla w_\para^\e\>\\
&\les \left|\RE \<\Lambda^s \nabla
  w_\para^\e,\Lambda^s(v_\para^\e\cdot \nabla^2 w_\para^\e)\>\right|+
  \left|\RE \<\Lambda^s \nabla 
  w_\para^\e,\Lambda^s(w_\para^\e\cdot \nabla^2 v^\e)\>\right|\\
&+ \rho^\e_\para(t)^2 + \rho^\e_\para(t)\left\|R_\para \nabla^2(\id
  -J_{1/\para})\Lambda^s\( |a_\para^\e|^2 - 
    |a^\e|^2\)\right\|_{L^2}\\
& + \rho^\e_\para(t)\left\|\nabla
  S_\para^\e\right\|_{H^s}. 
\end{align*}
Write the first term of the right hand side as:
\begin{align*}
  \RE \<\Lambda^s \nabla
  w_\para^\e,\Lambda^s(v_\para^\e\cdot \nabla^2 w_\para^\e)\>&= \RE
  \<\Lambda^s \nabla 
  w_\para^\e,v_\para^\e\cdot \nabla^2 \Lambda^s w_\para^\e\>\\
&+ \RE \<\Lambda^s \nabla
  w_\para^\e,\Lambda^s(v_\para^\e\cdot \nabla^2 w_\para^\e)-
  v_\para^\e\cdot \nabla^2 \Lambda^s w_\para^\e\>.
\end{align*}
Integration by parts and Kato-Ponce estimates \eqref{Kato-Ponce} yield:
\begin{align*}
  \left|\RE \<\Lambda^s \nabla
  w_\para^\e,\Lambda^s(v_\para^\e\cdot \nabla^2 w_\para^\e)\>\right|& \les
  \rho^\e_\para(t)^2+ \rho^\e_\para(t)\|\nabla
  v_\para^\e\|_{L^\infty}\|\nabla^2 
  w_\para^\e\|_{H^{s-1}}\\
& + \rho^\e_\para(t)\|\nabla v_\para^\e\|_{H^{s-1}}\|\nabla^2
  w_\para^\e\|_{L^\infty}\les  \rho^\e_\para(t)^2, 
\end{align*}
where we have used Corollary~\ref{cor:borneunif0} and Sobolev
embeddings. Similarly,
\begin{align*}
  \RE \<\Lambda^s \nabla
  w_\para^\e,\Lambda^s(w_\para^\e\cdot \nabla^2 v^\e)\>&= \RE
  \<\Lambda^s \nabla 
  w_\para^\e,w_\para^\e\cdot \nabla^2 \Lambda^s v^\e\>\\
&+ \RE \<\Lambda^s \nabla
  w_\para^\e,\Lambda^s(w_\para^\e\cdot \nabla^2 v^\e)-
  w_\para^\e\cdot \nabla^2 \Lambda^s v^\e\>,
\end{align*}
and:
\begin{align*}
  \left|\RE \<\Lambda^s \nabla
  w_\para^\e,\Lambda^s(w_\para^\e\cdot \nabla^2 v^\e)\>\right|& \les
  \rho^\e_\para(t)^2+ \rho^\e_\para(t)\|\nabla
  w_\para^\e\|_{L^\infty}\|\nabla^2 
  v^\e\|_{H^{s-1}}\\
& + \rho^\e_\para(t)\|\nabla w_\para^\e\|_{H^{s-1}}\|\nabla^2
  v^\e\|_{L^\infty}\les  \rho^\e_\para(t)^2.
\end{align*}
We also have
\begin{align*}
 \left\|R_\para \nabla^2(\id
  -J_{1/\para})\Lambda^s\( |a_\para^\e|^2 - 
    |a^\e|^2\)\right\|_{L^2}\les  \left\|\Lambda^s\( |a_\para^\e|^2 - 
    |a^\e|^2\)\right\|_{L^2}\les \rho^\e_\para(t),
\end{align*}
and we infer:
\begin{equation}\label{eq:restedw}
 \frac{d}{dt}\|\nabla w_\para^\e\|_{H^s}^2 \les \rho^\e_\para(t)^2 +
  \rho^\e_\para(t)\left\|\nabla 
  S_\para^\e\right\|_{H^s}. 
\end{equation}
Proceeding similarly for $d^\e_\para$, we find:
\begin{equation}\label{eq:rested}
 \frac{d}{dt}\|d_\para^\e\|_{H^s}^2 \les \rho^\e_\para(t)^2 +
  \rho^\e_\para(t)\left\|\Sigma_\para^\e\right\|_{H^s}. 
\end{equation}
Summing over \eqref{eq:resteLinf} and the time integrated
Equations~\eqref{eq:restedw}  and \eqref{eq:rested}, we complete the
first task of the program announced above:
\begin{equation}
  \label{eq:rgronwall}
  \rho^\e_\para(t)\les \int_0^t \rho^\e_\para(\tau)d\tau +
 \int_0^t\(\| S_\para^\e 
 (\tau)\|_{L^\infty}+\| \nabla S_\para^\e
 (\tau)\|_{H^s }+\| \Sigma_\para^\e
 (\tau)\|_{H^s}\)d\tau. 
\end{equation}
Since $S_\para^\e \in H^s$, Lemma~\ref{lem:sobolev} implies:
\begin{equation*}
  \rho^\e_\para(t)\les \int_0^t \rho^\e_\para(\tau)d\tau + \int_0^t\(\| \nabla
 S_\para^\e  (\tau)\|_{H^s }+\| \Sigma_\para^\e (\tau)\|_{H^s}\)d\tau. 
\end{equation*}
It is an easy consequence of Corollary~\ref{cor:borneunif0} and
Lemma~\ref{lem:cvop} that we have:
\begin{equation*}
  \| \nabla
 S_\para^\e\|_{L^\infty_T H^s }+\|
 \Sigma_\para^\e\|_{L^\infty_TH^s}\to 0 \quad\text{as }\para \to 0. 
\end{equation*}
We infer from Gronwall lemma that $\rho^\e_\para\to 0$ as $\para \to 0$,
uniformly on $[0,T]$. Therefore, we have:
\begin{equation*}
  (v^\e,a^\e)\in C([0,T];X^{s+1}\times X^{s})\quad ;\quad |a^\e|^2-c
  \in C([0,T];L^2),
\end{equation*}
and the existence part of Proposition~\ref{prop:EPoisson0} follows by
a bootstrap argument (to prove the extra smoothness). 
\smallbreak

Uniqueness follows  from the above computations: up to changing the
notations, we have the same estimates as above, with now
$S^\e=\Sigma^\e\equiv 0$. Uniqueness then follows from Gronwall
lemma. 

\smallbreak

To see that there exists $\phi^\e$ such that $v^\e = \nabla \phi^\e$,
apply the $\operatorname{curl}$ operator to the equation satisfied by
$v^\e$ \eqref{eq:systav0}. Energy estimates then show that
$  \nabla \times v^\e \equiv 0$.
We conclude thanks to \cite[Prop.~1.2.1]{JYClivre}.

\smallbreak

Before being more precise about the properties of $\phi^\e$ (we
already know that $\nabla \phi^\e \in C([0,T];X^\infty)$), we examine
the Poisson potential $\poi^\e$. We have
\begin{equation*}
  \Delta  \poi^\e = q\( |a^\e|^2-c\)\in C([0,T];H^\infty). 
\end{equation*}
We infer from Lemma~\ref{lem:infini} that 
\begin{equation*}
  {\mathcal F}_{y\to \xi}\(\nabla  \poi^\e\) \in C([0,T];L^1).
\end{equation*}
We deduce $\nabla  \poi^\e \in C([0,T]\times\R^n)$, and
Riemann-Lebesgue lemma implies that
\begin{equation*}
  \nabla  \poi^\e (t,x)\Tend {|x|}{\infty} 0. 
\end{equation*}
So far, we have worked with $\nabla  \poi^\e$ only, and we know that
it is smooth. At this stage, $\poi^\e$ is determined up to
a function of time only. The condition $\poi^\e (t,0)=0$ fixes the
value of that function, and yields a unique, smooth, Poisson potential
(so far, only its gradient was unique). 
As announced in the introduction, we explain why we cannot (in general)
impose the behavior
\begin{equation}\label{eq:poiinf}
  \poi^\e (t,x)\Tend {|x|}{\infty} 0,
\end{equation}
instead of $\poi^\e (t,0)=0$. We know from Lemma~\ref{lem:infini} that
for all $t\in 
[0,T]$, $\nabla \poi^\e (t,\cdot)\in L^p(\R^n)$, for $p>2n/(n-2)$. We can
then apply Lemma~\ref{lem:hormander} only when $n\ge 5$. The
following example shows that in space dimension $n=3$, we may have
$\nabla f(x)\Tend {|x|}{\infty} 0$, $\Delta f\in H^\infty$, and $f
(x)\Tend {|x|}{\infty} +\infty$:
\begin{equation*}
  f(x)=\log \( 1+|x|^2\),\quad x\in \R^3.
\end{equation*}
Note also that in the case $c\in L^1(\R^n)$ discussed below, we have
the additional property $\Delta \poi^\e \in C([0,T];L^1\cap
H^\infty)$, which makes it possible to impose \eqref{eq:poiinf}. 
Back to $\phi^\e$, we
have:
\begin{equation*}
  \nabla\( \d_t \phi^\e + \nabla \eik\cdot
  \nabla \phi^\e+\frac{1}{2}|\nabla \phi^\e|^2 +V_{\rm pert} +
  \poi^\e\)=0. 
\end{equation*}
We infer:
\begin{equation*}
   \d_t \phi^\e + \nabla \eik\cdot
  \nabla \phi^\e+\frac{1}{2}|\nabla \phi^\e|^2 +V_{\rm pert} +
  \poi^\e= F, 
\end{equation*}
where $F=F(t)$ is a function of time only.
In the above equation, all the terms are uniquely
determined, except $\d_t \phi^\e$ and $F$. 
Imposing $\phi^\e_{\mid t=0} =\phi_0$, and replacing $\phi^\e$ with
$\phi^\e +\int_0^t G(\tau)d\tau$ if necessary, we may assume that
$F\equiv 0$. This condition fully determines $\phi^\e$. 
This completes the
proof of Theorem~\ref{theo:existence}. 


\section{Convergence as $\e \to 0$: proof of
Theorem~\ref{theo:convergence}} 
\label{sec:convergence}
In this section, we prove
Theorem~\ref{theo:convergence}. First, the existence of  $(a,\phi)$
solving \eqref{eq:systlim1} follows from the proof of
Theorem~\ref{theo:existence}, since \eqref{eq:systlim1} is nothing but
\eqref{eq:systcomplet2} with $\e=0$. 
Denote 
  \begin{equation*}
   w_v^\e := v^\e -v=\nabla \phi^\e
    -\nabla \phi\quad ;\quad  w_a^\e := a^\e -a .
  \end{equation*}
The pair $(w_v^\e, w_a^\e)$ solves a system similar to
\eqref{eq:reste}: 
\begin{align*}
  &\d_t w_v^\e + w_v^\e\cdot \nabla v
  + (v_{\ei}+v^\e)\cdot \nabla w_v^\e  + \nabla\( \poi^\e - \poi \)=0.\\
&\d_t w_a^\e + w_v^\e\cdot \nabla a
  + (v_{\ei}+v^\e)\cdot \nabla  w_a^\e  
+ \frac{1}{2}\(w_a^\e\nabla\cdot  v^\e
  + a\nabla\cdot  w_v^\e\) 
 =i\frac{\e}{2}\Delta a^\e.\\
&\Delta \(\poi^\e - \poi \)=q\(|a^\e|^2 -|
  a|^2 \).\\
&\nabla\(\poi^\e - \poi \)(t,x)\to 0 \text{ as }|x|\to
  \infty\quad ;\quad w_v^\e\big|_{t=0} = 0\quad ; \quad
  w_a^\e\big|_{t=0} = r^\e. 
\end{align*}
Let $s>n/2+1$. Mimicking the computations made in
Section~\ref{sec:cvschemelin}, we find:
\begin{align*}
  \frac{d}{dt}\(\|w_v^\e(t)\|_{X^{s+1}}^2 +\|w_a^\e(t)\|_{H^s}^2  \)
  &\lesssim \|w_v^\e(t)\|_{X^{s+1}}^2 +\|w_a^\e(t)\|_{H^s}^2\\
+\e \|\Delta a^\e\|_{H^s}\|w_a^\e(t)\|_{H^s}
&+\|\nabla\(\poi^\e - \poi \)\|_{X^{s+1}}\|w_v^\e(t)\|_{X^{s+1}}.
\end{align*}
From Theorem~\ref{theo:existence},
\begin{equation*}
  \|\Delta a^\e\|_{L^\infty_TH^s}\lesssim 1.
\end{equation*}
To estimate the term corresponding to the Poisson potentials, write:
\begin{equation*}
  \|\nabla\(\poi^\e - \poi \)\|_{X^{s+1}} \les \|\nabla\(\poi^\e - \poi
  \)\|_{L^\infty} +  \|\Delta\(\poi^\e - \poi \)\|_{H^{s}}.
\end{equation*}
We estimate the first term thanks to Lemmas~\ref{lem:hormander} and
\ref{lem:sobolevp}: we have
\begin{equation*}
  \Delta \( \poi^\e - \poi \) = q\(|a^\e|^2-|a|^2\)=
  q\(|a^\e|^2-c+c-|a|^2\).  
\end{equation*}
Therefore, $\Delta \( \poi^\e - \poi \)$ is bounded in
$L^\infty_TH^s$. We see
from \eqref{eq:17h30} that 
$\d_{jk}^2 \( \poi^\e - \poi \)$ is bounded in
$L^\infty_TH^s$ for every pair $(j,k)$. Lemma~\ref{lem:hormander}
(with $p=2$) shows that there exists a function $\gamma^\e_j(t)$ of time
only such that
\begin{equation*}
  \d_{j} \( \poi^\e - \poi \)(t,\cdot)-\gamma^\e_j(t)\in
  L^{\frac{2n}{n-2}}(\R^n). 
\end{equation*}
On the other hand, for all $t\in [0,T]$, 
\begin{equation*}
  \d_{j} \( \poi^\e - \poi \)(t,x)\Tend {|x|}{\infty} 0.
\end{equation*}
Therefore, $\gamma^\e_j(t)\equiv 0$, and $\d_{j} \( \poi^\e - \poi
  \)(t,\cdot)\in L^{\frac{2n}{n-2}}(\R^n)$. The critical Sobolev
  embedding then shows that
  \begin{equation*}
    \|\nabla\(\poi^\e - \poi  \)\|_{L^{\frac{2n}{n-2}}}\les \|\Delta\(\poi^\e
    - \poi  \)\|_{L^2}.
  \end{equation*}
Along with Lemma~\ref{lem:sobolev} (with $q=2n/(n-2)$), this yields:
\begin{equation*}
  \|\nabla\(\poi^\e - \poi  \)\|_{L^\infty} \les \|\Delta\(\poi^\e
    - \poi  \)\|_{H^s},
\end{equation*}
and we have:
\begin{equation*}
  \|\nabla\(\poi^\e - \poi \)\|_{X^{s+1}} \les \|\Delta\(\poi^\e - \poi
  \)\|_{H^{s}}\les \|w_a^\e\|_{H^s}.  
\end{equation*}
We infer: 
\begin{align*}
  \frac{d}{dt}\(\|w_v^\e(t)\|_{X^{s+1}}^2 +\|w_a^\e(t)\|_{H^s}^2  \)
  \lesssim \|w_v^\e(t)\|_{X^{s+1}}^2 +\|w_a^\e(t)\|_{H^s}^2
 + \e \|w_a^\e(t)\|_{H^s}.
\end{align*} 
By assumption,
\begin{equation*}
  \|w_v^\e(0)\|_{X^{s+1}} +\|w_a^\e(0)\|_{H^s}  =\|r^\e\|_{H^s}\Tend
  \e 0 0,
\end{equation*}
and we conclude with Gronwall lemma:
\begin{equation*}
  \|w_v^\e\|_{L^\infty_TX^{s+1}} +\|w_a^\e\|_{L^\infty_TH^s}\les \e +
  \|r^\e\|_{H^s}. 
\end{equation*}
The strong convergence of the quadratic quantities described in
Theorem~\ref{theo:convergence} follows easily. Note that a
similar convergence has been obtained by P.~Zhang \cite{ZhangSIMA},
when $V_{\rm ext}\equiv 0=\alpha_0$ (hence $\eik\equiv
\beta_0$) and $c\in L^1(\R^n)$. The convergence in 
\cite{ZhangSIMA} is proved is a weaker sense though (in the sense
of measures), due to a different technical approach based on the use
of Wigner measures. 
\bigbreak

To conclude this section, we note that one must not expect $a
e^{i(\phi +\eik)/\e}$ to be a good \emph{pointwise} approximation of
$u^\e = a^\e e^{i(\phi^\e +\eik)/\e}$. We have:
\begin{align*}
  u^\e-ae^{i(\eik +\phi)/\e} &= a^\e e^{i(\eik +\phi^\e)/\e}-
  ae^{i(\eik +\phi)/\e} \\
& = \(a^\e-a\) e^{i(\eik +\phi^\e)/\e} +a e^{i\eik /\e}\(e^{i\phi^\e/\e}-
  e^{i\phi/\e}\).
\end{align*}
The first term is $\GO(\e)$ in $L^2\cap L^\infty$ (we avoid
differentiation because of rapid oscillations). The modulus of the last
term is of order
\begin{equation*}
  |a| \left|\sin \( \frac{\phi^\e -\phi}{\e}\)\right|.
\end{equation*}
Note that our results do not allow us to estimate the argument of the
sine function. Formally, it should not be 
smaller than $\GO(1)$ in general, so we must not expect $ae^{i(\eik
  +\phi)/\e}$ to be 
a good approximation for $u^\e$. To have a good approximation, we
would have to
compute the next term in the asymptotic expansion for $(a^\e,\phi^\e)$
as $\e\to 0$. We leave out this question at this stage here, because
we do not have completely satisfactory answers for that issue, and
resume this discussion when $c\in L^1(\R^n)$ below, a case where we
have more precise information at hand.


\part{Subquadratic eikonal phase}\label{part:2}


We now allow the external potential and the initial phase
to have quadratic components. After some geometrical reductions, the
analysis boils down to the previous one. This reveals some differences
though: for instance, even if $|a_0^\e|^2 - c\in L^2$, one must not expect
$|u^\e(t)|^2-c \in L^2$ for $t>0$. 
\begin{hyp}\label{hyp:general} Recall that $n\ge 3$.\\
\noindent$\bullet$  \emph{External potential:}  $\ext\in
C^\infty(\R\times\R^n)$ writes 
\begin{equation*}
  \ext (t,x)= V_{\rm quad}(t,x) + V_{\rm pert}(t,x),
\end{equation*}
where $V_{\rm quad}\in C^\infty(\R\times\R^n)$ is a polynomial of
degree at most two in $x$ ($\nabla^3 V_{\rm quad}\equiv 0$), and
$\nabla V_{\rm pert}\in C(\R;H^\infty)$.\\

\noindent$\bullet$  \emph{Doping profile:} it is a short range
perturbation of a constant. For simplicity, we assume that this
constant is $1$:
\begin{equation*}
  c = 1 + \widetilde c,\quad\text{where }\widetilde c \in
  L^1\cap H^\infty.
\end{equation*}    
 
\noindent$\bullet$ \emph{Initial amplitude:} it has the following
expansion,
\begin{equation*}
  a_0^\e(x) =a_0(x) + r^\e(x),
\end{equation*}
where $ a_0\in X^\infty$ is such that $|a_0|^2- 1\in
L^2(\R^n)$, and $r^\e\in H^\infty$, with 
\begin{equation*}
  \|r^\e\|_{H^s}\Tend \e 0 0,\quad \forall s\ge 0.
\end{equation*}

\noindent$\bullet$ \emph{Initial phase:} we have $\Phi_0\in
C^\infty(\R^n)$ with
\begin{equation*}
  \Phi_0(x) = \phi_{\rm quad}(x) +\phi_0(x),
\end{equation*}
where $\phi_{\rm quad}$ is a polynomial of order at most two, and
$\nabla \phi_0\in H^\infty$.
\end{hyp}
\begin{example}[External potential]
  We may take 
  \begin{equation*}
    V_{\rm quad}(t,x) = \sum_{j=1}^{n}\lambda_j(t) x_j^2,
  \end{equation*}
an anisotropic harmonic potential with smooth time-dependent coefficients. 
Of course, we may take $V_{\rm pert}\in C^\infty(\R;H^\infty)$. 
\end{example}

\section{The eikonal phase and the associated transport operator}
\label{sec:eik}

The generalization of Lemma~\ref{lem:eik0} is:
\begin{lemma}\label{lem:eik}
  Under the Assumption~\ref{hyp:general}, there exists $T^*>0$ and a unique
  solution $\eik\in C^\infty([0,T^*]\times\R^n)$ to:
  \begin{equation}
    \label{eq:eik}
    \d_t \eik +\frac{1}{2}|\nabla \eik|^2 +V_{\rm quad}(t,x)=0\quad ;\quad
    \phi_{{\rm eik} \mid t=0}=\phi_{\rm quad}\, .
  \end{equation}
This solution is a polynomial of order at most two in $x$:
$\nabla^3\eik \equiv 0$. 
\end{lemma}
\begin{proof}
  The first part of the lemma was established in
  \cite{CaBKW}. Consider 
  the Hamiltonian flow associated to 
  \begin{equation*}
    \frac{1}{2}|\xi|^2 + V_{\rm quad}(t,x),
  \end{equation*}
which yields $x(t,y)$ and $\xi(t,y)$ solving:
\begin{equation}
  \label{eq:hamilton}
\left\{
  \begin{aligned}
   &\d_t x(t,y) = \xi \(t,y\)\quad ;\quad x(0,y)=y,\\ 
   &\d_t \xi(t,y) = -\nabla_x V_{\rm quad}\(t,x(t,y)\)\quad ;\quad
   \xi(0,y)=\nabla \phi_{\rm quad}(y).
  \end{aligned}
\right.
\end{equation}
Following this flow and using a global inversion theorem (see
\cite{DG} for these general results), we construct $\eik$, locally in time,
but globally in space. The idea for the global inversion is to notice
that $\nabla_y x$ is the identity, plus a perturbation which is uniformly
bounded in space, and continuous in time with initial value equal to zero:
there exists $T^*>0$ such that, for all 
$t\in [0,T^*]$, $y\mapsto x(t,y)$ is a global diffeomorphism. We
denote by $y(t,x)$ its inverse. This
yields $\eik \in C^\infty([0,T^*]\times\R^n)$, with 
\begin{equation}\label{eq:nablaeik}
  \nabla \eik(t,x) =  \xi\(t,y(t,x)\).
\end{equation}
As a byproduct, the function $\eik$ is
sub-quadratic: $\d_x^\alpha \eik \in L^\infty([0,T^*]\times\R^n)$ as
soon as $|\alpha|\ge 2$.

Differentiating
\eqref{eq:eik} three times with respect to any triplet of space
variables, we see that $\Psi = \nabla^3 \eik$ solves a system of the
form:
\begin{equation*}
  \(\d_t +\nabla\eik\cdot \nabla\)\Psi = {\mathcal M}\Psi\quad ;\quad
  \Psi_{\mid t=0} =0,
\end{equation*}
where ${\mathcal M} \in L^\infty([0,T^*]\times \R^n)$ (${\mathcal M}$ is
a linear combination of derivatives of order at least two of
$\eik$). Note that the absence of source term and initial datum
follows from Assumption~\ref{hyp:general}. Since $\nabla \eik$ is
given by \eqref{eq:nablaeik}, we can then use the method
of characteristics: setting $\widetilde \Psi(t,y)= \Psi(t,x(t,y))$
(which makes sense since $x(t,\cdot)$ is a global diffeomorphism),
the above equation becomes
\begin{equation*}
  \d_t \widetilde \Psi = \widetilde {\mathcal M}\widetilde \Psi\quad ;\quad
  \widetilde \Psi_{\mid t=0} =0\quad ;\quad \widetilde{\mathcal M} \in
  L^\infty([0,T^*]\times \R^n). 
\end{equation*}
We conclude with Gronwall lemma that $\widetilde\Psi\equiv \Psi\equiv
0$.  

\smallbreak
Alternatively, one can prove that $\Psi\equiv 0$ by an elementary
integration by parts argument. Namely, since  
$\nabla^2\eik \in L^\infty([0,T^*]\times\R^n)$ and $\Psi(0,\cdot)=0\in
L^{2}(\R^{n})$, we have  
$\Psi\in C([0,T^*];L^{2}(\R^{n}))$ together with the energy
identity: for all $t\in [0,T^*]$,  
$$
\int\left\lvert \Psi(t,x)\right\rvert ^{2}\, dx
=\int_{0}^{t}\int \bigl(\Delta \eik (t',x) +2
\mathcal{M}(t',x)\bigr)\left\lvert \Psi(t',x)\right\rvert ^{2}\,
dxdt'. 
$$
Hence, again, the desired result follows from Gronwall lemma.
\end{proof}


In view of the energy estimates performed in Section~\ref{sec:lin}, we
will not consider \eqref{eq:eik}, but a nonlinear perturbation of this
equation. Indeed, if we try to mimic the computations after
Lemma~\ref{lem:contlinfini0}, and after having changed variables to
work on the characteristics, we have to estimate
$D_ta_\para^\e$ in $L^2$. From the equation,
\begin{align*}
  \|D_t a_\para^\e\|_{L^2}&\les \bigl(1+\left\lVert v_\para^\e
   \right\rVert_{L^\infty}+\left\lVert\nabla v_\para^\e
   \right\rVert_{L^2} \bigr)  
\bigl(\left\lVert a_\para^\e \right\rVert_{L^\infty}+
\left\lVert\nabla a_\para^\e \right\rVert_{L^2}\bigr)+
\left\lVert \Delta a_\para^\e  \right\rVert_{L^2}\\
&+ \|a_\para^\e \Delta
   \eik\|_{L^2}. 
\end{align*}
The last term is new, since now $\Delta
   \eik$ is a non-trivial function (of time only). This
   means that we must not even expect the last term to be finite! To
   overcome this difficulty, we proceed as on the baby model
   \begin{equation*}
     \d_t a + \frac{1}{2}a \Delta \eik =0,
   \end{equation*}
where from Lemma~\ref{lem:eik}, $\Delta \eik$ is a function of time only.
It is convenient to introduce the auxiliary function
\begin{equation*}
  \widetilde a(t,x) = a(t,x) \exp\(\frac{1}{2}\int_0^t \Delta
  \eik(\tau)d\tau\).  
\end{equation*}
Therefore, it is tempting to replace the condition $|u^\e|-1\in
L^\infty_TL^2$ with a condition of the form
\begin{equation*}
  \left| u^\e e^{\frac{1}{2}\int_0^t \Delta
  \eik(\tau)d\tau}\right|^2-1 \in L^\infty_TL^2.
\end{equation*}
Apparently, we have solved the issue mentioned above, but the price to
pay is that we no longer consider the quantity which is natural in
view of the Poisson equation. The idea is then to introduce a ``ghost
Poisson potential'':
\begin{equation*}
  \poi^\e = \widetilde \poi^\e + V_{\rm g}+ V_{\widetilde c}, \quad
  \text{where} 
\end{equation*}
\begin{align*}
  \Delta \widetilde \poi^\e &= qe^{-\int_0^t \Delta
  \eik(\tau)d\tau}\(\left| u^\e e^{\frac{1}{2}\int_0^t \Delta
  \eik(\tau)d\tau}\right|^2-1\),\\
\Delta V_{\rm g} &= q \(e^{-\int_0^t \Delta
  \eik(\tau)d\tau} -1\),\\ 
\Delta V_{\widetilde c} &= q \widetilde c.
\end{align*}
In particular, $\Delta V_{\rm g}$ is a function of time only: $V_{\rm
  g}$ is quadratic in $x$, and we may choose
\begin{equation}
  \label{eq:vg}
  V_{\rm g}(t,x) = q\frac{|x|^2}{2n}\(e^{-\int_0^t \Delta
  \eik(\tau)d\tau} -1\).
\end{equation}
Following the idea of \cite{CaBKW}, it is consistent to replace
$V_{\rm quad}$ with $V_{\rm quad}+V_{\rm g}$ in \eqref{eq:eik}, since
$V_{\rm g}$ is quadratic and cannot be considered as a perturbation or
a source term. Even though $V_{\rm g}$ depends on $\eik$, it is
reasonable to try to extend Lemma~\ref{lem:eik}. Indeed, if we
consider the iterative scheme
\begin{align*}
  \d_t \eik^{(j+1)} +\frac{1}{2}\left|\nabla \eik^{(j+1)}\right|^2 +V_{\rm
    quad}(t,x)&=-q\frac{|x|^2}{2n}\(e^{-\int_0^t \Delta
  \eik^{(j)}(\tau,x)d\tau} -1\),\\
     \phi^{j+1}_{{\rm eik} \mid t=0}&=\phi_{\rm quad},
\end{align*}
with $\eik^{(0)}= \phi_{\rm quad}$, we see that applying
Lemma~\ref{lem:eik} inductively shows that every iterate is a smooth,
sub-quadratic function. We have precisely:
\begin{proposition}\label{prop:eik}
  Under Assumption~\ref{hyp:general}, there exists $T^*>0$ and a unique
  solution $\eik\in C^\infty([0,T^*]\times\R^n)$, polynomial of order
  at most two in $x$
($\nabla^3\eik \equiv 0$),  to:
\begin{equation}
    \label{eq:eik2}
    \begin{aligned}
     \d_t \eik +\frac{1}{2}|\nabla \eik|^2 +V_{\rm
  quad}(t,x)&=-q\frac{|x|^2}{2n}\(e^{-\int_0^t \Delta 
  \eik(\tau)d\tau} -1\),\\
    \phi_{{\rm eik} \mid t=0}&=\phi_{\rm quad}\, .
    \end{aligned}
    \end{equation}
We denote:
\begin{equation}
  \label{eq:trace}
  \ga(t)\defn \frac{1}{2}\int_0^t \Delta  \eik(\tau)d\tau.
\end{equation}
\end{proposition}
\begin{proof}
  Inspired by Lemma~\ref{lem:eik}, we seek directly $\eik$ of the form
  \begin{equation*}
    \eik(t,x)= \,^{t}xM(t)x + \alpha(t)\cdot x + \beta(t), 
  \end{equation*}
where $M\in {\mathcal M}_{n\times n}(\R)$, $\alpha\in
\R^n$ and $\beta\in \R$. 
Plugging this expression into \eqref{eq:eik2} and identifying the
coefficients of the polynomials in $x$, we find:
\begin{align*}
  \dot M(t) + 2M(t)^2 + Q(t) &= -\frac{q}{2n}\(e^{-2\int_0^t {\rm
  Tr}M(\tau)d\tau} -1\)I_n\quad ;\quad M(0)=M_0,\\ 
\dot \alpha(t)+2M(t)\alpha(t)+E(t)&=0\quad ;\quad
  \alpha(0)=\alpha_0,\\
\dot \beta(t) +\frac{1}{2}|\alpha(t)|^2 +  \gamma(t)&=0\quad ;\quad
  \beta(0)=\beta_0,
\end{align*}
where
\begin{equation*}
  V_{\rm quad}(t,x)=\,^{t}xQ(t)x + E(t)\cdot x + \gamma(t)\quad ;\quad 
\phi_{\rm quad}(x) = \,^{t}xM_0 x +\alpha_0\cdot x+\beta_0.
\end{equation*}
Introducing the unknown function $R(t)=\int_0^t M(\tau)d\tau$, we see that
the equation in $M$ can be solved thanks to Cauchy-Lipschitz Theorem
applied to $(M(t),R(t))$. Then $\alpha(t)$ and $\beta(t)$ follow by
simple integration. 
\end{proof}
The above proof shows that unless $Q(t)\equiv 0=M_0$ (a case which
boils down to Part~\ref{part:1}), $\ga$ is a non-trivial function of
time. 
\smallbreak

The previous result implies that the characteristics associated to the
transport operator $\d_{t}+\nabla\eik\cdot\nabla$ present in
\eqref{eq:systcomplet2} can be described 
very easily.  
\begin{corollary}\label{cor:charac}
Let $x(t,y)$ be as defined in \eqref{eq:hamilton}. 
There exist $\alpha\in C^{\infty}([0,T^*];\R^n)$, and $A\in
C^{\infty}([0,T^*];{\mathcal M}_{n\times n}(\R))$, 
symmetric,  such that, for all 
$t\in [0,T^*]$,
$$
x(t,y)=e^{A(t)}y+\alpha(t).
$$
\end{corollary}
\begin{proof}
By Proposition~\ref{prop:eik}, $\nabla\eik(t,x)=M(t)x+\alpha(t)$ for some
symmetric matrix $M$. Since 
$$
\partial_{t}x(t,y)=\nabla\eik (t,x(t,y)),
$$
the result follows by integration. 
\end{proof}
\begin{remark}
Under Assumption~\ref{hyp:lin}, $\nabla\eik$ is a function
of time only, and the transport operator $\d_t+\nabla\eik\cdot \nabla$
is trivial. In the above proof, $M\equiv 0$, and we have $x(t,y)=
y+\int_0^t\alpha(\tau)d\tau$. This relation is reminiscent of
Avron--Herbst formula (see e.g. \cite{Cycon}). 
\end{remark}

\section{Main results}
\label{sec:road}
The analogue of Theorem~\ref{theo:existence} is:
\begin{theorem}\label{theo:const}
 Let Assumption~\ref{hyp:general} be  satisfied. 
 There exists
$T>0$ independent 
  of $\e \in ]0,1]$ and a  solution $u^\e \in
  L^\infty([0,T]\times \R^n)$
 to   \eqref{eq:schrod}-\eqref{eq:CI}, with
 \begin{equation*}
   \nabla\widetilde\poi^\e(t,x)= \nabla\( \poi^\e -V_{\rm
   g}-V_{\widetilde c}\)(t,x)\Tend {|x|}{\infty} 0,\quad 
   \widetilde\poi^\e(t,0)= 0,
 \end{equation*}
and such that $|u^\e e^\ga|^2-1 \in
  L^\infty([0,T];L^2)$, where $\ga$ is given by
  \eqref{eq:trace}. Moreover, one can write $u^\e = a^\e e^{i(\eik 
    +\phi^\e)/\e}$, where:
  \begin{itemize}
  \item $a^\e \in  C^\infty([0,T]\times\R^n)\cap C([0,T];X^\infty)$,
  and $|a^\e e^\ga|^2-1 \in 
  C([0,T];L^2)$. 
  \item $\eik$ is given by Proposition~\ref{prop:eik}.
  \item $\phi^\e \in C^\infty([0,T]\times\R^n)$ and $\nabla\phi^\e \in
  C([0,T];X^\infty)$. 
  \item We have the following uniform estimate: for every $s>n/2$,
  there exists $M_s$ \emph{independent of $\e \in ]0,1]$} such that
  \begin{equation*}
    \| a^\e\|_{L^\infty (0,T;X^s)}
    +\left\||a^\e e^\ga|^2-1\right\|_{L^\infty (0,T;L^2)}+ \|\nabla
    \phi^\e\|_{L^\infty (0,T;X^s)}\le M_s. 
  \end{equation*}
  \end{itemize} 
\end{theorem}
\begin{remark}
  We impose conditions on $\widetilde \poi^\e$, and not on
  $\poi^\e$. This is related to the \emph{arbitrary} choice \eqref{eq:vg} to
  integrate the ``ghost Poisson equation'' (this equation introduces
  additional degrees of freedom), since we will impose
  \begin{equation*}
    \nabla V_{\widetilde c}(x)\Tend {|x|}\infty 0\quad ;\quad
    V_{\widetilde c}(x)\Tend {|x|}\infty 0. 
  \end{equation*}
\end{remark}
Note that since $\ga$ is non-trivial, the above result shows that one
must not expect $|u^\e(t)|^2-1\in L^2$ for $t>0$. 
\smallbreak

Proceeding like before, we want $(a^\e,v^\e)$ to solve:
\begin{equation}\label{eq:systav}
\left\{
\begin{aligned}
 & \d_t v^\e + \(v_{\ei}+v^\e\)\cdot
  \nabla v^\e+v^\e\cdot \nabla v_{\ei} +\nabla V_{\rm pert} +
  \nabla\poi^\e=0,\\
&\d_t a^\e +\(v_{\ei}+v^\e\) \cdot \nabla a^\e +\frac{1}{2}a^\e
  \nabla\cdot \(v_{\ei}+v^\e\)
  = i\frac{\e}{2}\Delta a^\e,\\ 
 &\Delta \poi^\e =q\big(|a^\e|^2-c\big),
\end{aligned}
\right.
\end{equation}
together with
\begin{equation}\label{ci:systav}
\nabla\widetilde\poi^\e(t,x)\Tend {|x|}{\infty}0\quad ;\quad
  \widetilde\poi^\e(t,0)=0\quad ;\quad v^\e_{\mid 
  t=0}=\nabla\phi_0 \ ;\ a^\e_{\mid t=0}=a_0^\e . 
\end{equation}
With this existence result, we can study the asymptotic behavior as
$\e\to 0$ of the solution we construct:
\begin{theorem}\label{theo:convergence2}
  Under Assumption~\ref{hyp:general}, there
  exists a smooth solution $(a,\phi)$ of \eqref{eq:systav} with
  $\e=0$,  
such that 
  $a,\nabla\phi \in C([0,T],X^\infty)$, $|a e^\ga|^2-1 \in
  C([0,T],L^2)$, 
  and 
  $$
  \|a^\e-a\|_{L^\infty_{T}H^s}+ \|\nabla(\phi^\e
  -\phi)\|_{L^\infty_{T}X^s}\Tend \e 0 0, \qquad \forall s>n/2.
  $$
  In particular, the position density and the momentum density
  converge: 
  \begin{align*}
&|u^\e|^2 \Tend \e 0 |a|^2\quad \text{in } L^\infty_TH^s, \text{ and}\\
&\e\IM\(\overline{u^\e
     e^{-i\eik/\e}}  \nabla\(u^\e   e^{-i\eik/\e}\)\) \Tend \e 0
     |a|^2\nabla\phi \quad \text{in }L^\infty_TX^s,\ \forall s>n/2.
  \end{align*}
\end{theorem}
\begin{remark}
 We slightly altered the usual notion of momentum density, by removing
 first the eikonal phase $\eik$. Indeed, we do not prove that
 \begin{equation*}
   |a^\e|^2 \nabla \eik \Tend \e 0 |a|^2 \nabla \eik \quad \text{in
    }L^\infty_TX^s,
 \end{equation*}
since $\nabla \eik$ may grow linearly in $x$, while $|a|^2$ morally
goes to $1$ as $|x|\to \infty$, hence
$|a|^2 \nabla \eik  \not \in L^\infty_TX^s$. 
\end{remark}

We first show  
that the solutions of \eqref{eq:systav} exist and
are uniformly bounded for a time interval independent of $\e$. 
\begin{proposition}\label{prop:EPoisson} 
Let Assumption~\ref{hyp:general} be satisfied. Let $s>n/2$. 
For all $ M>\ M_{0}>0$, there exists $0<T\le T^*$
such that, if for all $\e\in 
[0,1]$,
\begin{equation}\label{TM:CI}
\left\lVert \nabla
  \phi_{0}\right\rVert_{H^{s+2}}+ 
\bigl\lVert \left\lvert a_{0}^{\e} \right\rvert
^2-1\bigr\rVert_{L^{2}}+\left\lVert 
  a_{0}^{\e}\right\rVert_{X^{s+1}}\le  M_{0}, 
\end{equation}
then the Cauchy problem~\eqref{eq:systav} has a
unique classical solution   
$(v^{\e},a^{\e})$ in $ C^{\infty}([0,T]\times\R^{n})$ such that
\begin{equation}\label{TH:NORM}
\left\lVert  v^{\e}\right\rVert_{L^\infty_T X^{s+2}}
+ \left\lVert \left\lvert a^{\e} e^{\ga}\right\rvert
^2-1\right\rVert_{L^\infty_T L^{2}} 
+ \left\lVert  a^{\e}\right\rVert_{L^\infty_TX^{s+1}}\le M.
\end{equation}
\end{proposition}

We perform some geometrical reductions so that the proofs of the above
results follow from Section~\ref{sec:lin}. 

\section{Reduction to the first case}\label{sec:reduction}
We begin by proving that \eqref{eq:systav} is 
equivalent to a system which does not involve the operator
$v_{\ei}\cdot\nabla$, thanks to Corollary~\ref{cor:charac}. Resuming the
notations of Section~\ref{sec:eik}, define, for any function $f$ of
time and space:
\begin{equation*}
  \widetilde f(t,y) = f\(t,x(t,y)\). 
\end{equation*}
Working with $\widetilde f$ instead of $f$, the characteristics
associated to $v_{\ei}\cdot\nabla$ are straightened so that:
\begin{equation*}
  \d_t \widetilde{f}(t,y)= (\d_t+v_{\ei}\cdot\nabla
)f\(t,x(t,y)\).
\end{equation*}
The good news for us is the fact that the above change of
variable does not change the structural
properties of \eqref{eq:systav}.  
Indeed, Corollary~\ref{cor:charac} implies that 
\begin{equation}\label{eq:charac}
\widetilde{\nabla f}(t,y) = e^{-A(t)}\nabla \widetilde{f}(t,y),  
\end{equation}
for some symmetric $n\times n$ matrix $A(t)$ which is independent of $y$.


We are now in position to make precise the fact that the change of
variables does not change the  structural properties. 
\begin{lemma}\label{lem:structure}
Fix $t\in [0,T^*]$ and set $\delta_t \defn e^{-A(t)}\nabla$, where $A$
is as in Corollary~\ref{cor:charac}.  
The following properties hold:

\noindent $(1)$ For all $u\in H^{2}(\R^n)$ and all $v\in
W^{1,\infty}(\R^n)$ one has: 
\begin{align*}
\RE\scal{i\delta^*_t\delta_t u}{u}=0\quad ;\quad 
2\scal{v\cdot\delta_t u}{u}=-\scal{(\delta_t v)u}{u}. 
\end{align*}
\noindent $(2)$ The Fourier multiplier $\nabla
(-\delta^*_t\delta_t)^{-1}\delta_t$ is 
well defined and bounded on Sobolev spaces:  
for all $\sigma\ge 0$, there exists a constant $K_\sigma$ independent of
$t\in [0,T^*]$ such that:
\begin{equation*}
\left\lVert \nabla (-\delta^*_t\delta_t)^{-1}\delta_t
  u\right\rVert_{H^\sigma}\le K_\sigma \left\lVert u\right\rVert_{H^\sigma},
\qquad \forall u\in H^{\sigma}(\R^n). 
\end{equation*}
\noindent $(3)$ For all function $u\colon\R^n\rightarrow \R$,
\begin{equation*}
u(x)\Tend{|x|}{\infty}0 \quad\Leftrightarrow\quad
u(x(t,y))\Tend{|y|}{\infty}0. 
\end{equation*}
\end{lemma}
\begin{proof}
By integrating by parts, the first property follows from the fact
that $\delta_t$ is a linear combination of spatial derivatives  
whose coefficients are constant symmetric matrices. The property $(2)$
is immediate using Fourier transform.  The property $(3)$ is obvious. 
\end{proof}

\begin{notation}
Introduce the operator $\d$ by, for all $u\colon
[0,T^*]\rightarrow\mathcal{S}'(\R^n)$, 
\begin{equation*}
(\d u)(t)\defn e^{-A(t)}\nabla u(t).
\end{equation*}
\end{notation}
The difference between the above notation and
Lemma~\ref{lem:structure} is that $\delta_t$ is defined for \emph{fixed}
$t\in [0,T^*]$. Following what we did in Section~\ref{sec:eik},
introduce
\begin{equation*}
  \widetilde a^\e \defn a^\e e^\ga\quad ;\quad \poi^\e =\widetilde
  \poi^\e +V_{\rm g}+V_{\widetilde c}. 
\end{equation*}
Since $\widetilde c\in L^1\cap H^\infty$, $\Delta^{-1}\widetilde c$ is well
defined as a temperate distribution:
\begin{equation*}
  \Delta^{-1}\widetilde c = -\F^{-1}\(|\xi|^{-2}\F( \widetilde c) \). 
\end{equation*}
Setting $\widetilde V_{\rm pert} \defn V_{\rm pert} +q
\Delta^{-1}\widetilde c$, 
we still have $\nabla \widetilde V_{\rm pert}\in C(\R;H^\infty)$, from
Lemma~\ref{lem:poisson}. 
With these notations, \eqref{eq:systav} is equivalent to:
\begin{equation}\label{eq:systfinal}
\left\{
\begin{aligned}
 & \d_t v^\e + v^\e\cdot\d v^\e+
 v^\e\cdot \d v_{\ei} 
 +\d \widetilde V_{\rm pert} +
  \d\widetilde \poi^\e=0,\\
&\d_{t} \widetilde a^{\e}+{v^\e}\cdot\d
 \widetilde a^\e+\frac{1}{2}\widetilde a^\e\d\cdot v^{\e}
= -i\frac{\e}{2}\d^{*}\d \widetilde a^\e,\\
& \d^{*}\d \widetilde \poi^\e =-qe^{-2\ga}\bigl( \lvert \widetilde a^\e\rvert
 ^2-1\bigr).  
\end{aligned}
\right.
\end{equation}
Note that the fact that the right hand side of the equation for $a^\e$
is skew-symmetric remains, from the first point of
Lemma~\ref{lem:structure}: following the idea of E.~Grenier 
\cite{Grenier98}, this is crucial in the proof of
Theorem~\ref{theo:existence}. This is so thanks to
Corollary~\ref{cor:charac}, and would not be if $\eik$ was not exactly
polynomial. 

\smallbreak

Directly from Corollary~\ref{cor:charac}, we verify that, for all
$\sigma\ge 0$, there exists $C_\sigma$ such that, for all $t\in [0,T^*]$, 
\begin{equation}\label{normequiv}
C_\sigma^{-1}\left\lVert u(x(t,\cdot)) \right\rVert_{H^\sigma} \le
\left\lVert u\right\rVert_{H^\sigma}\le C_\sigma \left\lVert
  u(x(t,\cdot)) \right\rVert_{H^\sigma} , \quad 
\forall u\in H^{\sigma}(\R^n).
\end{equation}
Similarly, there exists a constant $C$ such that
\begin{equation*}
C^{-1}\left\lVert u(x(t,\cdot)) \right\rVert_{L^\infty} \le
\left\lVert u\right\rVert_{L^\infty}\le C\left\lVert
  u(x(t,\cdot)\right\rVert_{L^\infty} , \quad 
\forall u\in L^{\infty}(\R^n).\notag
\end{equation*}
Note that \eqref{eq:systfinal} is very similar to
\eqref{eq:systav0}. The transport operator is simplified, but we have
two new features:
\begin{itemize}
\item The term $v^\e\cdot \d v_{\ei} $ in the equation for $v^\e$.
\item The factor $e^{-2\ga}$ in the Poisson equation. 
\end{itemize}
The latter changes very little computations, since $\ga$ is a
function of time only. One can check that the term $v^\e\cdot \d
v_{\ei}$ does not require a modification of the proof given in
Section~\ref{sec:lin}; unlike for $\d_t a^\e$, we do not estimate $\d_t
v^\e$ in $L^2$. Therefore, Theorem~\ref{theo:const}
follows. Similarly, Theorem~\ref{theo:convergence2} follows like in
Section~\ref{sec:convergence}, up to some slight modifications; like
for the proof of Theorem~\ref{theo:const}, replace
$(a^\e,\poi^\e,V_{\rm pert})$
with $(\widetilde a^\e,\widetilde \poi^\e, \widetilde V_{\rm pert})$.


\part{Integrable doping profile}\label{part:3}


\section{Integrable doping profile}\label{sec:integrable}

There are many results concerning the case when the doping profile $c$
is decaying at spatial infinity, say  
$c\in L^1(\R^n)$. We refer for instance to \cite{ZhangSIMA} and
references therein. We restrict our
attention to the case $n=3$ for simplicity. 
\begin{hyp}\label{hyp:int} We consider the case $n= 3$.\\
\noindent$\bullet$  \emph{External potential:}  $\ext\in
C^\infty(\R\times\R^3)$ writes 
\begin{equation*}
  \ext (t,x)= V_{\rm quad}(t,x) + V_{\rm pert}(t,x),
\end{equation*}
where $V_{\rm quad}\in C^\infty(\R\times\R^3)$ is a polynomial of
degree at most two in $x$ ($\nabla^3 V_{\rm quad}\equiv 0$), and
$V_{\rm pert}\in C(\R;H^\infty)$.\\

\noindent$\bullet$  \emph{Doping profile:}  $c\in L^1(\R^3)\cap 
X^\infty$.   \\ 
 
\noindent$\bullet$ \emph{Initial amplitude:} $ a_0^\e(x) =a_0(x) + \e
a_1(x)+\e r_1^\e(x)$, 
where $ a_0,a_1,r_1^\e\in H^\infty$, with 
\begin{equation*}
  \|r_1^\e\|_{H^s}\Tend \e 0 0,\quad \forall s\ge 0.
\end{equation*}

\noindent$\bullet$ \emph{Initial phase:} we have $\Phi_0\in
C^\infty(\R^n)$ with
\begin{equation*}
  \Phi_0(x) = \phi_{\rm quad}(x) +\phi_0(x),
\end{equation*}
where $\phi_{\rm quad}$ is a polynomial of order at most two, and
$\phi_0\in X^\infty$. 
\end{hyp}

Our goal is to state a convergence result which is more precise
than Theorem~\ref{theo:convergence}: we shall need some properties of
$\phi^\e$, and not only $\nabla \phi^\e$. This is why we change the
boundary conditions to solve the Poisson equation: we consider
\begin{equation}
  \label{eq:SPint}
\begin{aligned}
      i\e \d_t u^\e +\frac{\e^2}{2}\Delta u^\e &= \ext u^\e +
      \poi^\e u^\e ,\quad (t,x)\in \R\times \R^3,\\
\Delta \poi^\e = q\(|u^\e|^2-c\)&,\quad \nabla\poi^\e(t,x)\to 0\text{
      and }\poi^\e(t,x)\to 0 \text{ as
      }|x|\to \infty,\\
u^\e_{\mid t=0} &= a_0^\e(x)e^{i\Phi_0(x)/\e}.
\end{aligned} 
\end{equation}
With these boundary conditions (which are as in
\cite{CastoM3AS97,ZhangSIMA} for instance), we can define $\Delta^{-1}$ as:
\begin{equation}\label{eq:fond}
  \Delta^{-1}f = -\frac{1}{4\pi |x|}\ast f. 
\end{equation}

\begin{theorem}\label{theo:existence2}
Let $n=3$. Under Assumption~\ref{hyp:general}, assume furthermore that $V_{\rm
  pert}\in C(\R;X^\infty)$, $c\in L^{1}(\R^3)$, $a_0\in L^2(\R^3)$ and
  $\phi_0\in  L^\infty(\R^3)$. There exists 
$0<T\le T^*$ independent 
  of $\e \in ]0,1]$ and a unique solution $u^\e \in
  C^\infty([0,T]\times \R^3)\cap
 C([0,T];H^\infty)$  
 to   \eqref{eq:SPint}. Moreover, one can write $u^\e =
 a^\e e^{i(\eik 
    +\phi^\e)/\e}$, where:
  \begin{itemize}
  \item $a^\e\in  C([0,T];H^\infty)$. 
  \item $\eik$ is given by Lemma~\ref{lem:eik}. 
  \item $\phi^\e \in C([0,T];X^\infty)$. 
  \item We have the following uniform estimate: for every $s\ge 0$,
  there exists $M_s$ \emph{independent of $\e \in ]0,1]$} such that
  \begin{equation*}
    \| a^\e\|_{L^\infty_{T}H^{s}}  +\|
    \phi^\e\|_{L^\infty_{T}X^s}\le M_s. 
  \end{equation*}
  \end{itemize}
\end{theorem}
Note that existence and
uniqueness for \eqref{eq:SPint} can be established in a larger class
of functions, thanks 
to Strichartz estimates. We refer for instance to \cite{CastoM3AS97}
for the case 
with no external potential, and simply recall that similar Strichartz
estimates are available in the presence of a smooth, subquadratic
external potential (\cite{Fujiwara}, see also \cite{CaCCM}). Note also
that the term $\Delta^{-1}c$ can be treated as a ``nice'' linear
potential, thanks to Lemma~\ref{lem:poisson} and the following:
\begin{lemma}\label{lem:fond}
  The operator $\Delta^{-1}$ defined by \eqref{eq:fond} maps
  $L^1\cap L^2(\R^3)$ to $\F (L^1(\R^3))$, where $\F$ denotes the
  Fourier transform. Moreover, there exists $C$ such that
  \begin{equation*}
    \left\|\F \( \Delta^{-1}\varphi\)\right\|_{L^1}\le
    C\(\|\varphi\|_{L^1} + \|\varphi\|_{L^2} \),\quad \forall \varphi
    \in L^1\cap L^2(\R^3). 
  \end{equation*}
\end{lemma}
Uniqueness for \eqref{eq:SPint} follows easily:
\begin{align*}
  i\e \d_t u^\e +\frac{\e^2}{2}\Delta u^\e = \ext u^\e +
      q\Delta^{-1}\(|u^\e|^2-c\) u^\e.
\end{align*}
Let $u^\e$ and $v^\e$ be two solutions in $C([0,T^\e];H^\infty)$ of
the above equation, with the same initial data, for some
$T^\e>0$. Note that the dependence upon $\e$ is irrelevant, since
$\e>0$ is fixed. The difference $w^\e = u^\e-v^\e$ solves
\begin{align*}
  i\e \d_t w^\e +\frac{\e^2}{2}\Delta w^\e = \ext w^\e +
     q\Delta^{-1}\(|u^\e|^2-c\) w^\e+
     q\Delta^{-1}\(|u^\e|^2-|v^\e|^2\) v^\e. 
\end{align*}
The basic energy estimate yields:
\begin{align*}
  \e \frac{d}{dt}\|w^\e\|_{L^2}^2 &\les
  \left\|\Delta^{-1}\(|u^\e|^2-|v^\e|^2\) v^\e \right\|_{L^2}\|w^\e\|_{L^2} \\
&\les \left\|\Delta^{-1}\(|u^\e|^2-|v^\e|^2\)\right\|_{L^\infty} \|v^\e
  \|_{L^2}\|w^\e\|_{L^2}\\
&\les \( \left\|  |u^\e|^2-|v^\e|^2\right\|_{L^1} + \left\|
  |u^\e|^2-|v^\e|^2\right\|_{L^2}\) \|v^\e
  \|_{L^2}\|w^\e\|_{L^2}\\
&\les \( \|u^\e\|_{L^2} +\|v^\e\|_{L^2} +\|u^\e\|_{L^\infty}
  +\|v^\e\|_{L^\infty} \)\|w^\e\|_{L^2}^2 \|v^\e \|_{L^2}.  
\end{align*}
Uniqueness then follows from the Gronwall lemma.
\smallbreak

To prove the existence part of Theorem~\ref{theo:existence2}, we
consider
 \begin{equation}
  \label{eq:systcomplet2int}
\begin{aligned}
  \d_t \phi^\e + \nabla \eik\cdot
  \nabla \phi^\e+\frac{1}{2}|\nabla \phi^\e|^2 +V_{\rm pert} + \poi^\e&=0
  \ ;\ \phi^\e_{\mid t=0}=\phi_0.\\
\d_t a^\e +\nabla \(\phi^\e+\eik\) \cdot \nabla a^\e +\frac{1}{2}a^\e \Delta
  \(\phi^\e+\eik\) &= i\frac{\e}{2}\Delta a^\e\ ;\ a^\e_{\mid
  t=0}=a_0^\e.\\ 
\Delta \poi^\e =q\big(|a^\e|^2-c\big)\quad ;\quad \nabla\poi^\e(t,x)\Tend
  {|x|}{\infty}0&\text{ and }\poi^\e(t,x)\Tend
  {|x|}{\infty}0. 
\end{aligned}
\end{equation}
The geometrical reduction presented in Section~\ref{sec:reduction}
makes it possible to transform the transport operator $\d_t
+v_{\ei}\cdot\nabla$ into $\d_t$. Unlike in
Section~\ref{sec:reduction}, we may keep  the term $\Delta
\eik$. Since $\Delta \eik$ is a function of time only, and since we
work with $a^\e\in C([0,T];H^s)$, the term $a^\e \Delta \eik$ can be
treated like a perturbative term. 
\smallbreak

Since the proof of Theorem~\ref{theo:existence2} involves more
classical arguments, we essentially skip it, so that we can focus our
discussion on the semi-classical limit $\e \to 0$. 

After the geometrical reduction, \eqref{eq:systcomplet2int} becomes
what we would have found directly in the case $V_{\rm
  quad}=0=\phi_{\rm quad}$, up to terms which can be treated by
Gronwall lemma. We may for instance resume the approach of
Section~\ref{sec:lin}, and replace $X^s$ with $H^s$. This way, we
construct $a^\e, v^\e\in C([0,T];H^\infty)$.

To complete the proof of Theorem~\ref{theo:existence2}, we finally
notice that $\phi^\e \in C([0,T];L^\infty)$, from Lemma~\ref{lem:fond}
and \eqref{eq:systcomplet2int} integrated along the characteristics.
\smallbreak

We can  now establish the analogue of Theorem~\ref{theo:convergence},
with a pointwise description. To do so, we introduce the solution to
\begin{equation}
  \label{eq:systlim1int}
 \begin{aligned}
  \d_t \phi + \nabla \eik\cdot
  \nabla \phi+\frac{1}{2}|\nabla \phi|^2 +V_{\rm pert} + \poi&=0
  \ ;\ \phi_{\mid t=0}=\phi_0.\\
\d_t a +\nabla \(\phi+\eik\) \cdot \nabla a +\frac{1}{2}a \Delta
  \(\phi+\eik\) &= 0\ ;\ a_{\mid
  t=0}=a_0.\\ 
\Delta \poi =q\big(|a|^2-c\big)\quad ;\quad \nabla\poi(t,x)\Tend
  {|x|}{\infty}0&\text{ and }\poi(t,x)\Tend
  {|x|}{\infty}0.
\end{aligned} 
\end{equation}
This system has a unique solution $(\phi,a)\in C([0,T];X^\infty\times
H^\infty)$. As pointed out at the end of
Section~\ref{sec:convergence}, the triplet $(\eik,\phi,a)$ does not
suffice to describe the pointwise limit of $u^\e$ as $\e \to
0$. This is the reason why in Assumption~\ref{hyp:int}, we want to
know $a^\e$ up to $o(\e)$ instead of $o(1)$ only. 
Consider the linearized system:
\begin{equation}
  \label{eq:systlim2}
 \begin{aligned}
  \d_t \phi_1 + \nabla (\eik+\phi)\cdot
  \nabla \phi_1 + V&=0
  \ ;\ \phi_{1\mid t=0}=0.\\
\d_t b +\nabla \(\eik+\phi\) \cdot \nabla b +\frac{1}{2}b \Delta
  \(\eik+\phi\)+\\
+\nabla \phi_1\cdot \nabla a + \frac{1}{2}a \Delta
  \phi_1 = \frac{i}{2}&\Delta a\ ;\ b_{\mid t=0}=a_1.\\ 
\Delta V =2q\operatorname{Re}\(\overline{a}b\)\  ;\ \nabla V(t,x)\Tend
  {|x|}{\infty}0&,\text{ and }V(t,x)\Tend
  {|x|}{\infty}0. 
\end{aligned} 
\end{equation}
It has a unique solution $(\phi_1,b)\in  C([0,T];X^\infty\times
H^\infty)$. 
\begin{theorem}\label{theo:7.4}
  Under the Assumption~\ref{hyp:int}, the solution to \eqref{eq:SPint} can
  be approximated at leading order by
  $ae^{i\phi_1}e^{i(\eik+\phi)/\e}$:
  \begin{equation*}
    \left\|u^\e -
    ae^{i\phi_1}e^{i(\eik+\phi)/\e}\right\|_{L^\infty_T(L^2\cap 
    L^\infty)} \to 0\text{ as }\e \to 0.
  \end{equation*}
\end{theorem}
\begin{remark}
  In general, $\phi_1$ is not trivial provided that $a_1\not
  \equiv 0$, and the amplitude of $u^\e$ is, at leading order,
  $ae^{i\phi_1}$. This phenomenon is due to the fact that from the
  point of view of geometric optics, \eqref{eq:schrod}-\eqref{eq:CI}
  (or \eqref{eq:SPint})
  is supercritical: to describe the exact solution at leading order as
  in Theorem~\ref{theo:7.4}, it is necessary to know its initial data
  up to $o(\e)$. This phenomenon 
may lead to instability results as in
  \cite{CaARMA}: modifying $a^\e_0$ at order $\sqrt\e$ for instance,
  affects the solution $u^\e$ at order $\O(1)$ for times of order
  $\sqrt\e$.     
\end{remark}
\begin{proof}[Sketch of the proof]
  The idea is to resume the approach of
  Section~\ref{sec:convergence}. Set 
  \begin{equation*}
    \widetilde w_v^\e = \nabla \(\phi^\e-\phi -\e \phi_1\)\quad ;
    \quad \widetilde w_a^\e= a^\e-a -\e b.  
  \end{equation*}
Proceeding as in Section~\ref{sec:convergence}, we find, for $s$
sufficiently large:
\begin{equation*}
  \|\widetilde w_v^\e\|_{L^\infty_TX^{s+1}} +\|\widetilde
  w_a^\e\|_{L^\infty_TH^s}\les \e^2 + 
  \e\|r_1^\e\|_{H^s}. 
\end{equation*}
As above, we infer an $L^2$ estimate for $\widetilde w_v^\e$:
\begin{equation*}
  \|\widetilde w_v^\e\|_{L^\infty_TH^{s+1}} +\|\widetilde
  w_a^\e\|_{L^\infty_TH^s}\les \e^2 + 
  \e\|r_1^\e\|_{H^s}, 
\end{equation*}
and directly from the equation,
\begin{equation*}
 \| \phi^\e-\phi -\e \phi_1\|_{L^\infty([0,T]\times \R^3)}\les \e^2 + 
  \e\|r_1^\e\|_{H^s}=o(\e).
\end{equation*}
We conclude:
\begin{align*}
  \left| u^\e - ae^{i\phi_1}e^{i(\eik+\phi)/\e}\right|&=\left|
  a^\e e^{i\phi^\e/\e}  - ae^{i\phi_1}e^{i\phi/\e}\right| = \left|
  a^\e e^{i\phi^\e/\e}  - ae^{i(\phi+\e \phi_1)/\e}\right|\\
& \les \e |b|+   |\widetilde
  w_a^\e| + |a|\left| \sin\(\frac{\phi^\e-\phi -\e
  \phi_1}{\e}\)\right|\\
& \les \e |b|+   |\widetilde
  w_a^\e| + |a|\times o(1).
\end{align*}
The result follows by taking the $L^2$ or the $L^\infty$ norm in space.
\end{proof}

\bibliographystyle{amsplain}
\bibliography{poisson}

\end{document}